\documentclass[12pt,reqno]{amsart}
\usepackage{amssymb,delarray}
\usepackage{amsfonts}
\usepackage{epsfig}
\usepackage[all]{xy}

\def\leq{\leqslant}
\def\geq{\geqslant}

\newtheorem{thm}{Theorem}
\newtheorem{lem}
{Lemma}
\newtheorem{prop}
{Proposition}
\newtheorem{claim}
{Claim}

\newtheorem{cor}
{Corollary}
{Remark}
\newtheorem{que}
{Question}
\newtheorem{ex}{Example}

{\catcode`\@=11
\gdef\n@te#1#2{\leavevmode\vadjust{%
 {\setbox\z@\hbox to\z@{\strut#1}%
  \setbox\z@\hbox{\raise\dp\strutbox\box\z@}\ht\z@=\z@\dp\z@=\z@%
  #2\box\z@}}}
\gdef\leftnote#1{\n@te{\hss#1\quad}{}}
\gdef\rightnote#1{\n@te{\quad\kern-\leftskip#1\hss}{\moveright\hsize}}
\gdef\?{\FN@\qumark}
\gdef\qumark{\ifx\next"\DN@"##1"{\leftnote{\rm##1}}\else
 \DN@{\leftnote{\rm??}}\fi{\rm??}\next@}}

\begin{document}
\baselineskip=13.7pt plus 2pt 

\title[Factorizations in finite groups] {Factorizations in finite groups}
\author[Vik.S. Kulikov]{Vik.S. Kulikov}

\address{Steklov Mathematical Institute}
 \email{kulikov@mi.ras.ru}
\dedicatory{} \subjclass{}
\thanks{This
research was partially supported by grants of NSh-4713.2010.1, RFBR
11-01-00185, and by AG Laboratory HSE, RF government grant, ag.
11.G34.31.0023. }
\begin{abstract}
A necessary condition for uniqueness of factorizations of elements
of a finite group $G$ with factors belonging to a union of some
conjugacy classes of $G$ is given. This condition is sufficient if
the number of factors belonging to each conjugacy class is big
enough. The result is applied to the problem on the number of
irreducible components of the Hurwitz space of degree $d$ marked
coverings of $\mathbb P^1$ with given Galois group $G$ and fixed
collection of local monodromies.
\end{abstract}

\maketitle
\setcounter{tocdepth}{1}


\def\st{{\sf st}}



\section*{Introduction} \label{introduc}
Let $f:X\to \mathbb P^1$ be a morphism of a non-singular irreducible
projective curve $X$ (defined over the field of complex numbers
$\mathbb C$) onto the projective line $\mathbb P^1$. Denote by
$\mathbb C(X)$ the field of rational functions on $X$. The morphism
$f$ defines a finite extension $f^*: \mathbb C(z)\hookrightarrow
\mathbb C(X)$ of the field of rational functions $\mathbb C(\mathbb
P^1)\simeq \mathbb C(z)$. Denote by $G$ the Galois group of this
extension.

Let us choose a point $z_0\in \mathbb P^1$ such that $z_0$ is not a
branch point of $f$ and number the points of $f^{-1}(z_0)=\{
w_1,\dots ,w_d\}$, where $d=\deg f$. We will call the morphism $f$
with fixed numbering of the points of $f^{-1}(z_0)$ a {\it marked
covering}.

Let $z_1, \dots, z_n\in\mathbb P^1$ be the set of branch points of
$f$. The numbering of the points of  $f^{-1}(z_0)$ defines a
homomorphism $f_*:\pi_1(\mathbb P^1\setminus \{ z_1,\dots ,z_n\},
z_0)\to \Sigma_d$ of the fundamental group $\pi_1=\pi_1(\mathbb
P^1\setminus \{ z_1,\dots ,z_n\}, z_0)$ to the symmetric group
$\Sigma_d$. The image $\text{im} f_*$ acts transitively on
$f^{-1}(z_0)$ and it is isomorphic to $G$ (so we can identify
$\text{im} f_*$ and $G$). Let $\gamma_1,\dots,\gamma_n$ be simple
loops around, respectively, the points $z_1,\dots, z_n$ starting at
$z_0$ and such that they generate the group $\pi_1$. The image
$g_j=f_*(\gamma_j)\in G$ is called a {\it local monodromy} of $f$ at
the point $z_j$. Note that the set $\{ g_1,\dots, g_n\}$ of local
monodromies generate the group $G$. The local monodromy $g_j$
depends on the choice of $\gamma_j$, therefore it is defined
uniquely up to conjugation in $G$. Denote by $O=C_1\bigsqcup\dots
\bigsqcup C_m\subset G$ the union of conjugacy classes of all local
monodromies and by $\tau_i$ the number of local monodromies of $f$
belonging to the conjugacy class $C_i$. The pair $(G,O)$ is called
an {\it equipped group} and the collection
$\tau=(\tau_1C_1,\dots,\tau_mC_m)$ is called the {\it monodromy
type} of $f$.

Let $\text{HUR}_{d,G,\tau}^m(\mathbb P^1)$ be the Hurwitz space (see
the definition of Hurwitz spaces in \cite{F}) of marked degree $d$
coverings of $\mathbb P^1$ with Galois group $G$ and monodromy type
$\tau$. The famous Clebsch -- Hurwitz Theorem (\cite{Cl}, \cite{H})
states that if $G=\Sigma_d$ and $O$ is the set of transpositions,
then $\text{HUR}_{d,\Sigma_d,\tau}^m(\mathbb P^1)$ consists of a
single irreducible component if $\tau= (nO)$ with even $n\geq
2(d-1)$ and it is empty overwise. Generalizations of Clebsch --
Hurwitz Theorem were obtained in \cite{Ku1} and \cite{Ku2}. In
particular, in \cite{Ku2} it was proved that for an equipped group
$(\Sigma_d,O)$ with $O=C_1\sqcup\dots \sqcup C_m$, where $C_1$ is
the conjugacy class of an odd permutation leaving fixed at least two
elements, the Hurwitz space $\text{HUR}_{d,\Sigma_d,\tau}^m(\mathbb
P^1)$ is irreducible if $\tau_1$ is big enough. On the other hand,
the example in \cite{W} shows that
$\text{HUR}_{8,\Sigma_8,\tau}^m(\mathbb P^1)$ consists at least of
two irreducible components, where $\tau=(1C_1,1C_2,1C_3)$ and $C_1$
is the conjugacy class of permutation $(1,2)(3,4,5)$, $C_2$ is the
conjugacy class of $(1,2,3)(4,5,6,7)$, and $C_3$ is the conjugacy
class of $(1,2,3,4,5,6,7)$. Therefore we can not wait that for a
fixed equipped finite group $(G,O)$ the number of irreducible
components of  $\text{HUR}_{d,G,\tau}^m(\mathbb P^1)$ does not
depend on the monodromy type $\tau$. But, we can wait that this
number does not depend on $\tau$ if $\tau_i$ is big enough for some
$i$ such that the elements of $C_i$ generate the group $G$.

In subsection \ref{amb}, for each equipped finite group $(G,O)$ such
that the elements of $O$ generate $G$, we define a number
$a_{(G,O)}$ depending on $G$ and $O$, call it the {\it ambiguity
index} of $(G,O)$, and, as a straightforward corollary of Theorems
\ref{Cuniq1} and \ref{Cuniq2} (see subsection \ref{TT}) and results
of \cite{Ku1}, we have
\begin{thm} For each equipped finite group $(G,O)$, $O=C_1\sqcup\dots
\sqcup C_m$, such that the elements of $O$ generate the group $G$,
there is a constant $T$ such that the number of irreducible
components of each non-empty Hurwitz space
$\text{HUR}_{d,G,\tau}^m(\mathbb P^1)$ is equal to $a_{(G,O)}$ if
$\tau_i\geq T$ for all $i=1,\dots,m$.

If the elements of $O_1=C_1\sqcup\dots \sqcup C_k$ for some $k<m$
generate the group $G$, then there is a constant $T_1$ such that the
number of irreducible components of $\text{HUR}_{d,G,\tau}^m(\mathbb
P^1)$ is not more than $a_{(G,O_1)}$ if $\tau_i\geq T_1$ for
$i=1,\dots, k$.
\end{thm}

This article is a continuation of \cite{Ku1}, in which the
investigation of the factorization semigroups over finite groups was
started. For a convenience of the reader, the main definitions and
useful statements from \cite{Ku1} are reminded in section $1$. In
section $2$, to each equipped group $(G,O)$, we associate a
$C$-group whose factorization semigroup is the same as the
factorization semigroup of $(G,O)$ and we investigate a connection
between this $C$-group and $(G,O)$. In section $3$, we prove the
stability of the factorization semigroups over finite $C$-groups.
Theorem \ref{main1}, proved in this section, plays the key role in
the proof of all main results of this article (see section $4$). In
section 5, we give a solution of the word problem for finite
$C$-groups and give an algorithm of computation of the ambiguity
index $a_{(G,O)}$ for an equipped finite group $(G,O)$.

The author is grateful to O.V. Kulikova who draws attention of the
author on the important role of van Kampen Lemma in geometric group
theory, only due to the Lemma it was possible to prove Theorem
\ref{main1}.

\section{Semigroups over groups}
\subsection{Factorization semigroups} \label{semigr}
A pair $(G,O)$, where $G$ is a group and $O$ is a subset of $G$
invariant under the inner automorphisms, is called an {\it equipped
group}. In the sequel, we will assume that ${\bf{1}}\not\in O$ and
$O$ consists of a finite number of conjugacy classes $C_i$ of $G$,
$O= C_1\sqcup \dots\sqcup C_m$, and the numbering of these conjugacy
classes is fixed.

A homomorphism $f: G_1\to G_2$ is a {\it homomorphism of equipped
groups} $(G_1,O_1)$ and $(G_2,O_2)$ if $f(O_1)\subset O_2$.

A semigroup $S(G,O)$, generated by the letters of the alphabet
$X=X_O=\{ x_g \mid g\in O\}$ being subject to the relations
\begin{equation} \label{rel1} x_{g_1}\cdot
x_{g_2}=x_{g_2}\cdot\,x_{g_2^{-1}g_1g_2}=x_{g_1g_2g_1^{-1}}\cdot
x_{g_1},\end{equation} is called the {\it factorization semigroup
with factors in} $O$ (or a {\it factorization semigroup over the
group} $G$). A homomorphism $\alpha =\alpha_G:S(G,O)\to G$, given by
$\alpha (x_g)=g$ for each $x_g\in X$, is called the {\it product
homomorphism}. We will denote by $g_1\dots g_n$ the image
$\alpha(x_{g_1}\cdot \, .\, .\, .\,\cdot x_{g_n})$ of
$s=x_{g_1}\cdot \, .\, .\, .\,\cdot x_{g_n}\in S(G,O)$.

The action $\rho$ of the group $G$ on the set $X$, given by
$$ x_a\in X\mapsto \rho (g)(x_a)=x_{gag^{-1}}\in X,$$
 defines a homomorphism
$\rho : G\to \text{Aut} (S(G,O))$ . The action $ \rho (g)$ on
$S(G,O)$ is called the {\it simultaneous conjugation } by $g\in G$.
Put $\lambda (g)=\rho(g^{-1})$ and $\lambda _{S} =\lambda
\circ\alpha$, $\rho _{S} =\rho \circ\alpha$.

\begin{claim} {\rm (\cite{Ku1})} \label{cl1}
For all $s_1,\, s_2\in S(G,O)$ we have
\[
s_1\cdot s_2=s_2\cdot \lambda_S(s_2)(s_1)=\rho_S(s_1)(s_2)\cdot s_1.
\]
\end{claim}

To each element $s=x_{g_1}\cdot .\, .\, .\, \cdot x_{g_n}\in
S(G,O)$,  let us associate a positive integer $ln(s)=n$ called the
{\it length} of $s$. It is easy to see that $ln: S(G,O)\to \mathbb
Z_{\geq 0}=\{ {\bf{a}}\in \mathbb Z\mid {\bf{a}}\geq 0\}$ is a
homomorphism of semigroups.

For each $s=x_{g_1}\cdot\, .\, . \, . \, \cdot x_{g_n}\in S(G,O)$
denote by $G_s$ the subgroup of $G$ generated by the images $\alpha
(x_{g_1})=g_1,\dots ,\alpha (x_{g_n})=g_n$ of the factors $x_{g_1},
\dots , x_{g_n}$, and denote by $G_O$ the subgroup of $G$ generated
by the elements of $O$.

\begin{claim} {\rm (\cite{Ku1})} \label{cc} The subgroup $G_s$ of $G$ is well defined, that is,
it does not depend on a presentation of
$s$ as a product of  generators $x_{g_i}\in X_O$.
\end{claim}

For subgroups $H$ and $Z$ of a group $G$ denote by
$$S(G,O)^H=\{s\in S(G,O)\mid G_s=H\},$$
$$S(G,O)_{Z}= \{ s\in S(G,O)\mid \alpha (s)\in Z \},$$ and
$S(G,O)_{Z}^H= S(G,O)_{Z}\cap S(G,O)^H$. It is easy to see that
$S(G,O)^H$ (respectively $S(G,O)_{Z}^H$) is isomorphic to the
semigroup $S(H,H\cap O)^H$ (respectively, isomorphic to $S(H,H\cap
O)_{Z}^H$) and the isomorphism is induced by imbedding $(H,H\cap
O)\hookrightarrow (G,O)$.

\begin{prop} {\rm (\cite{Ku1})} \label{simple} Let $(G,O)$ be an equipped group and let $s\in S(G,O)$. We have
\begin{itemize} \item[($1$)] $\text{ker}\, \rho $ coincides with the centralizer
$C_{O}$ of the group $G_{O}$ in $G$;
\item[($2$)] if
$\alpha (s)$ belongs to the center $Z(G_s)$ of $G_s$, then for each
$g\in G_s$ the action $\rho (g)$ leaves fixed the element $s\in
S(G,O)$;
\item[($3$)] if $\alpha (s\cdot x_g)$ belongs to the center $Z(G_{s\cdot x_g})$ of
$G_{s\cdot x_g}$, then $s\cdot x_g=x_g\cdot s$,
\item[($4$)] if $\alpha (s)=\bf{1}$, then $s\cdot s'=s'\cdot s$
for any $s'\in S(G,O)$.
\end{itemize}
\end{prop}

\begin{claim}{\rm (\cite{Ku1})} \label{commut} For any equipped group $(G,O)$ the
semigroup $S(G,O)_{\bf{1}}$ is contained in the center of the
semigroup $S(G,O)$ and, in particular, it is a commutative
subsemigroup.
\end{claim}

It is easy to see that if $g\in O$ is an element of order $n$, then
$x_g^n\in S(G,O)_{\bf{1}}$.
\begin{lem} {\rm (\cite{Ku1})} \label{fried}  Let $s\in
S(G,O)_{Z(G_O)}$ and $s_1\in S(G,O)$ be such that $G_{s_1}=G_O$,
where $Z(G_O)$ is the center of $G_O$. Then
\begin{equation} \label{oo} s\cdot s_1= \rho(g)(s)\cdot s_1
\end{equation}
for all $g\in G_O$.

In particular, if  $s\in S(G,O)$ is such that $G_s=G$ and $C\subset
O$ is  a conjugacy class of $G$, then for any $g_1,g_2\in C$ we have
\begin{equation} \label{pp} x_{g_1}^{n}\cdot s=x_{g_2}^{n}\cdot
s\end{equation} if $g_1^n$ belongs to the center $Z(G)$ of $G$.
\end{lem}
\begin{prop} {\rm (\cite{Ku1})} \label{simpl}
The elements of $S(G,O)_{\bf{1}}^G$ are fixed under the conjugation
action of $G$.
\end{prop}

\section{$C$-groups and $C$-graphs}
\subsection{ $C$-graphs of $C$-groups} \label{C-rel}
Remind the definition of $C$-groups. By definition (see, for
example, \cite{Ku}), a $C$-group $G$ is an equipped group $(G,O)$
such that the elements of $O$ are the generators (so called,
$C$-generators) of the group $G$ being subject to the relations
\begin{equation} \label{C-rel} g_i^{-1}g_jg_i=g_k,\quad (g_i,g_j,g_k)\in M,\end{equation}
where $M$ is a
subset of $O^3$. A homomorphism $f:G_1\to G_2$ of $C$-groups is
called a $C$-homomorphism if it is a homomorphism of equipped
groups. In particular, two $C$-groups $G_1$ and $G_2$ are $C$-{\it
isomorphic} if they are isomorphic as equipped groups.

To each $C$-group, we can associate a directed graph, called
$C$-{\it graph.} To give a definition of $C$-graphs, consider a
directed graph $\Gamma=(V,E)$, where $V$ is the set of vertices of
$\Gamma$ and the set of its edges $E$ is a collection $\{
e_{v_i,v_j}=(v_i,v_j)\}$ of ordered pairs of its vertices (some of
edges can be loops, that is, the equality $v_i=v_j$ is allowed). For
each vertex $v\in V$, let us denote by $T_v=\{ e_{v,v_i}\}$ (resp.,
$H_v=\{ e_{v_i,v}\}$) the set of edges whose tails (resp., heads)
are the vertex $v$. A directed graph $\Gamma$ is called a $C$-{\it
graph} if each edge $e\in E$ is labeled by an element of $V$, that
is, there is a map $f:E\to V$ so that the label of $e\in E$ is
$f(e)\in V$ (in the sequel, an edge $e_{v_1,v_2}$ with label
$f(e_{v_1,v_2})=v$ will be denoted by $e_{v_1,v_2,v}$), and $\Gamma
$ is such that the following five conditions are satisfied:
\begin{itemize}
\item[($i$)] for each vertex $v\in V$ the restrictions $f_{\mid T_v}$, $f_{\mid H_v}$ of the map
$f$ to $T_v$ and $H_v$ are one to one correspondences with $V$;
\item[($ii$)] for each vertex $v$ the head $v_1$ of the edge
$e_{v,v_1,v}$ is the vertex $v$, that is, $e_{v,v_1,v}$ is the loop
($v_1=v$);
\end{itemize}
Note that, by condition $(i)$, the tail (resp., the head) $v_1$ and
the label $v_2$ define uniquely the edge $e$ whose tail (resp.,
head) is $v_1$ and whose label is $v_2$. Therefore a sequence $v_1,
\dots, v_n$ and a tail $v_0$ define uniquely a path
$l(v_0;v_1,\dots,v_n)$ starting at the vertex $v_0$ along edges (in
the  positive direction) with labels $v_1, \dots, v_n$.

The third and fourth conditions are
\begin{itemize}
\item[($iii$)] if for some two vertices $v_1$ and $v_2$ the edge
$e_{v_1,v_3,v_2}$ is a loop, that is, $v_1=v_3$, then the edge
$e_{v_2,v_4,v_1}$  is also a loop ($v_2=v_4$).
\item[($iv$)] for any
edge $e_{v_1,v_2,v_3}$ and for any vertex $v$ the ends of the paths
$l(v;v_1,v_3)$ and $l(v;v_3,v_2)$ coincide.
\end{itemize}

To each labeled directed graph $\Gamma$, let us associate a
two-dimensional complex $K_{\Gamma}$ whose $1$-skeleton is $\Gamma$
and whose two-sells are the quadrangles $Q_{(v,e_{v_1,v_2,v_3})}$
one to one corresponding to the pairs $(v,e_{v_1,v_2,v_3})\in
V\times E$ such that the border  $\partial Q_{(v,e_{v_1,v_2,v_3})}$
of $Q_{(v,e_{v_1,v_2,v_3})}$ is the loop $l(v;v_1,v_3)\cdot
l(v;v_3,v_2)^{-1}$.

The fifth  condition is
\begin{itemize}
\item[($v$)] if for some two edges
$e_{v,v_1,v_2}$ and $e_{v,v_1,v_3}$ the loop $e_{v,v_1,v_2}\cdot
e_{v,v_1,v_3}^{-1}$ represents the unity of the fundamental group
$\pi_1(K_{\Gamma},v)$, then $v_2=v_3$.
\end{itemize}

To each $C$-group $(G,O)$, one can associate a $C$-graph. By
definition, the $C$-{\it graph $\Gamma=\Gamma_{(G,O)}$ of a
$C$-group} $(G,O)$ is a $C$-graph whose set of vertices $V=\{
v_{g_i}\mid g_i\in O\}$ is in one to one correspondence with the set
$O$. Two vertices $v_{g_1}$ and $v_{g_2}$, $g_1, g_2\in O$ are
connected by a labeled edge $e_{v_{g_1},v_{g_2},v_g}$ if and only if
we have the relation $g^{-1}g_1g=g_2$ for some $g\in O$.

Conversely, to each $C$-graph $\Gamma$ one can associate a $C$-group
$G_{\Gamma}=(G,Y)$ the set of $C$-generators $Y=\{ y_{v_i} \mid
v_i\in V\}$ of $G_{\Gamma}$ is in one to one correspondence with the
set $V$ of vertices of $\Gamma$. In $G_{\Gamma}$ there is a relation
$y_{v_3}^{-1}y_{v_1}y_{v_3}=y_{v_2}$ if and only if there is an edge
$e_{v_1,v_2,v_3}$ in $\Gamma$.

\begin{claim} \label{iso4} For each $C$-graph $\Gamma$, the $C$-group
$G_{\Gamma}$ is $C$-isomorphic to $G_{\Gamma_{G_{\Gamma}}}$. For
each $C$-group $G$, the $C$-graphs $\Gamma_G$ and
$\Gamma_{G_{\Gamma_G}}$ are isomorphic. \end{claim} \proof
Obvious. \qed \\

In the sequel, for a $C$-group $G_{\Gamma}$ the generators
$x_{y_v}$, $v\in \Gamma$, of the semigroup $S(G_{\Gamma},Y)$ will be
denoted by $x_v$.

We say that a subgraph $\Gamma_1$  of a $C$-graph $\Gamma$ is a
$C$-{\it subgraph} if $\Gamma_1$  is a $C$-graph.

Let $\Gamma_1$ be a $C$-subgraph of a $C$-graph $\Gamma$. Consider
the $C$-groups $(G,Y)=G_{\Gamma}$ and $(G_1,Y_1)=G_{\Gamma_1}$ and
their factorization semigroups $S(G,Y)$ and $S(G_1,Y_1)$. The
embedding $i: \Gamma_1\hookrightarrow \Gamma$ defines  the natural
homomorphism $i_*:G_{\Gamma_1}\to G_{\Gamma}$ of $C$-groups and the
homomorphism $i_*: S(G_1,Y_1)\to S(G,Y)$ of their factorization
semigroups given, respectively, by $i_*(y_{v_i})=y_{v_i}$ and
$i_*(x_{v_i})=x_{v_i}$ for $v_i\in \Gamma_1\hookrightarrow\Gamma$.
\begin{claim} \label{is1}
Let $\Gamma_1$ be a $C$-subgraph of a $C$-graph $\Gamma$. Then the
homomorphism  $i_*: S(G_{\Gamma_1},Y_1)\to S(G_{\Gamma},Y)$ is an
embedding.
\end{claim}
\proof Obvious. \qed \\

For each $C$-graph $\Gamma$ there is a homomorphism from the
$C$-group $G_{\Gamma}$ to the automorphism group of $\Gamma$. The
action of $G_{\Gamma}$ is defined as follows: if $v^{\prime}\in
\Gamma$ then the action of the $C$-generator $y_{v^{\prime}}$ is
given by the rule: for a vertex $v$ of $\Gamma$ the image
$y_{v^{\prime}}(v)$ is the head of the edge with tail $v$ and label
$v^{\prime}$, for an edge $e_{v_1,v_2,v_3}$ the image
$y_{v^{\prime}}(e_{v_1,v_2,v_3})=e_{y_{v^{\prime}}(v_1),y_{v^{\prime}}(v_2),y_{v^{\prime}}(v_3)}$.
It follows from conditions ($i$) -- ($v$) of the definition of
$C$-graphs that this action is well defined.

In the sequel, we will consider only finitely generated $C$-groups
(as groups without equipment) and $C$-graphs consisting of finitely
many connected components. Denote by $m$ the number of connected
components of $C$-graph $\Gamma$. Then it is easy to see that
$G_{\Gamma}/[G_{\Gamma},G_{\Gamma}]\simeq \mathbb Z^m$ and any two
$C$-generators $y_{v_1}$ and $y_{v_2}$ are conjugated in the
$C$-group $G_{\Gamma}$ if and only if $v_1$ and $v_2$ belong to the
same connected component of $\Gamma$, that is, the set $Y$ of
$C$-generators of the $C$-group $G_{\Gamma}$ is the union of $m$
conjugacy classes of $G_{\Gamma}$. Denote by $\text{ab}:
G_{\Gamma}\to H_1(G_{\Gamma},\mathbb
Z)=G_{\Gamma}/[G_{\Gamma},G_{\Gamma}]$ the natural epimorphism. In
the sequel, we will assume that some numbering of the connected
components of $\Gamma$ is fixed. In this case the group
$H_1(G_{\Gamma},\mathbb Z)\simeq \mathbb Z^m$ has the natural base
consisting of vectors $\text{ab}(y_{v})=(0,\dots,0,1,0\dots,0)$,
where $1$ stands on the $i$-th place if $v$ belongs to the $i$th
connected component of $\Gamma$. Denote the composition $ab\circ
\alpha_{G_{\Gamma}}$ by $\tau$.

Let $l=l(v_0;v_1,\dots,v_n)$ be a path in a $C$-graph $\Gamma$. The
number $n$ is called the {\it length} of $l$. The smallest positive
integer $p_v$ (maybe, $n_v=\infty$), such that for any vertex $v_i$
of $\Gamma$ the path $l(v_i;v,\dots, v)$ of length $p_v$ is a loop
with origin and end at $v_i$, is called the {\it period} of $v$. It
is easy to see that $p_v=\min \{p\in \mathbb N\mid y_{v}^{p}\in
Z(G_{\Gamma})\}$, where $Z(G_{\Gamma})$ is the center of
$G_{\Gamma}$.

\begin{claim} \label{per} If $v_1$ and $v_2$ belong to the same
connected component of a $C$-graph $\Gamma$, then $p_{v_1}=p_{v_2}$.
\end{claim}
\proof The elements $y_{v_1}$ and $y_{v_2}$ are conjugated in
$G_{\Gamma}$. Therefore, if $y_{v_1}^{p_{v_1}}\in Z(G_{\Gamma})$
then $y_{v_2}^{p_{v_1}}=y_{v_1}^{p_{v_1}}$ and hence $y_{v_2}^{p_{v_1}}\in Z(G_{\Gamma})$. \qed \\

\begin{claim} \label{ppp} If a vertex $v_1$ of a $C$-graph $\Gamma$ is such that its period $p_{v_1}=1$,
then the $C$-group $G_{\Gamma}$ is naturally isomorphic to the
direct product $G_{\Gamma_1}\times \mathbb F_1$, where $\mathbb F_1$
is a free group generated by $y_{v_1}$ and the $C$-group
$G_{\Gamma_1}$ is generated by the all $C$-generators $y_{v}$, where
$v\neq v_1$, and it is associated with $C$-graph $\Gamma_1$ obtained
from $\Gamma$ deleting the vertex $v_1$ and the all edges labeled by
$v_1$.
\end{claim}
\proof Evident. \qed \\

A $C$-group $G_{\Gamma}$ is called a {\it finite $C$-group} if the
$C$-graph $\Gamma$ is a finite graph (note that a finite $C$-group
is infinite in the sense of usual groups).

Let $\Gamma=\Gamma_1\sqcup\dots\sqcup\Gamma_m$ be the decomposition
into the disjoint union of the connected components of a finite
$C$-graph $\Gamma$ and let $p_i$ be the period of the set of
vertices $V_i=\{ v_{i,1},\dots, v_{i,n_i}\}$ of the connected
component $\Gamma_i$. The element
$$c=\prod_{i=1}^m\prod_{j=1}^{n_i}y_{v_{i,j}}^{p_i}$$
is called the {\it canonical element} of the $C$-group $G_{\Gamma}$.

\begin{prop} \label{fincom} Let $G_{\Gamma}$ be a finite $C$-group. Then
the commutator $[G_{\Gamma},G_{\Gamma}]$ is a finite group.
Moreover, each element $g\in [G_{\Gamma},G_{\Gamma}]$ can be written
in the form
\begin{equation}\label{ccom} g=c^{-1}\prod_{i=1}^{m}y_{i,1}^{k_ip_i}y_{i,1}^{a_{i,1}}y_{i,2}^{a_{i,2}}\dots
y_{i,n_i}^{a_{i,n_i}},
\end{equation}
 where $c$ is the canonical element of $G_{\Gamma}$ and the integers $k_i$ and $a_{i,j}$
satisfy the following relations and inequalities
\begin{equation} \label{xxxx}
\sum_{j=1}^{n_i}a_{i,j}+k_ip_i=n_ip_i,\end{equation}
\begin{equation}\label{yyyy}
 0<a_{i,j}\leq p_i,  \qquad 0\leq k_i< n_i.
\end{equation}
\end{prop}
\proof Applying relations (\ref{C-rel}) and since
$y_{i,j}^{p_i}=y_{i,1}^{p_i}$, each element $g\in [G,G]$ can be
written in the form
\begin{equation}\label{xyzw} g=\prod_{i=1}^{m}y_{i,1}^{k_ip_i}y_{i,1}^{b_{i,1}}y_{i,2}^{b_{i,2}}\dots
y_{i,n_i}^{b_{i,n_i}},\end{equation} where the integers $k_i$ and
$b_{i,j}$ satisfy the following relations and inequalities
\begin{equation} \label{xxxxx}
\sum_{j=1}^{n_i}b_{i,j}+k_ip_i=0,\end{equation}
\begin{equation}\label{yyyyy}
 |b_{i,j}|\leq p_i-1,  \qquad |k_i|< n_i.
\end{equation}
For fixed integers $p_i$ and $n_i$ the set of integer solutions of
equations (\ref{xxxxx}) under restrictions (\ref{yyyyy}) is
finite. Therefore $[G_{\Gamma},G_{\Gamma}]$ is a finite group.
To obtain presentation (\ref{ccom}) from (\ref{xyzw}), it suffices
to multiply presentation (\ref{xyzw}) of $g$ by $c^{-1}c$ and one more to use relations $y_{i,j}^{p_i}=y_{i,1}^{p_i}$. \qed \\

\subsection{Canonical elements of factorization semigroups}
In notations used above, the element
$$\displaystyle s_{\Gamma^{\prime}}=\prod_{l=1}^{k}\prod_{j=1}^{n_{i_l}}
x_{v_{i_l,j}}^{p_{i_l}}\in S(G,Y),$$ is well defined and it is
called the {\it canonical element} associated with a subgraph
$\Gamma^{\prime}=\Gamma_{i_1}\sqcup \dots \sqcup\Gamma{i_k}$ of the
$C$-graph $\Gamma$. If $\Gamma^{\prime}=\Gamma$, then the element
$s_{\Gamma}$ is called the  {\it canonical element} of the semigroup
$S(G,Y)$. Obviously, $s_{\Gamma^{\prime}}$ belongs to the center of
$S(G,Y)$, since each its factor $x_{v_{i_l,j}}^{p_{i_l}}$ belongs to
the center of $S(G,Y)$.

An element $s_1\in S(G,Y)$ is said to be a {\it divisor} of an
element $s\in S(G,Y)$ if there is $s_2\in S(G,Y)$ such that
$s=s_1\cdot s_2$.
\begin{lem} \label{divis}  An element $s\in S(G,Y)$ of length
$ln(s)\leq k$ is a divisor of $s_{\Gamma^{\prime}}^k$ if $s$ can be
represented as a word in generators $x_{v_{i,j}}$, where $v_{i.j}\in
\Gamma^{\prime}$.
\end{lem}
\proof Obvious.
\begin{lem}\label{x1} Let $G_{\Gamma}=(G,Y)$ be a finite $C$-group.
Let the $i$th coordinate $\tau_i(s)$ of $\tau(s)$ for
an element $s=s^{\prime}\cdot s^{\prime\prime}\in S(G,Y)^G$ is not
less than $n_ip_i+1$. Assume that $\tau_i(s^{\prime\prime})=0$. Then
for $1\leq j\leq n_i$ the element $s$ can be written in the form:
$s=x_{i,j}^{p_i}\cdot s_i\cdot s^{\prime\prime}$, where $s_i\cdot
s^{\prime\prime}\in S(G,Y)^G$.
\end{lem}
\proof Since $\tau_i(s)\geq n_ip_i+1$, there are at least $p_i+1$
factors $x_{v_{i,j}}$ in a factorization of $s^{\prime}$,
$s^{\prime}=x_{v_{i_1,j_1}}\cdot \, .\, .\, . \cdot x_{v_{i_k,j_k}}$
having the same $i$ and $j$. Applying relations (\ref{rel1}), we can
move them to the left and after that we obtain a new factorization
$s^{\prime}=x_{v_{i,j}}^{p_i}\cdot (x_{v_{i,j}}\cdot s_i^{\prime})$.
Obviously, $s_i\cdot s^{\prime\prime}:=(x_{v_{i,j}}\cdot
s_i^{\prime})\cdot s^{\prime\prime}$ belongs to $S(G,Y)^G$, since
$s\in S(G,Y)^G$. Applying Lemma \ref{fried}, we complete the proof.
\qed
\begin{cor} \label{can}
Let $G_{\Gamma}=(G,Y)$ be a finite $C$-group. Let, for an element
$s=s^{\prime}\cdot s^{\prime\prime}\in S(G,Y)^G$, the $i$th
coordinate $\tau_i(s)$ is not less than $2n_ip_i+1$. Assume that
$\tau_i(s^{\prime\prime})=0$. Then the element $s$ can be written in
the form: $s=s_{\Gamma_i}\cdot s_i\cdot s^{\prime\prime}$, where
$s_i\cdot s^{\prime\prime}\in S(G,Y)^G$.
\end{cor}

\subsection{Ample subgraphs of $C$-graphs}
Let $\Gamma^{\prime}$ be a union of some connected components of a
$C$-graph $\Gamma$. We say that $\Gamma^{\prime}$ (and respectively,
the union of conjugacy classes of $C$-generators of the $C$-group
$G_{\Gamma}$ corresponding to $\Gamma^{\prime}$) is {\it ample} if
any two vertices of $\Gamma$ can be connected by a path along edges
of $\Gamma$ labeled by vertices belonging to $\Gamma^{\prime}$. In
language of $C$-groups, it means that any two conjugated
$C$-generators of $G_{\Gamma}$ are conjugated by some element of the
subgroup $G_{\Gamma^{\prime}}$ of $G_{\Gamma}$ generated by the
$C$-generators $y_v$, $v\in \Gamma^{\prime}$. Note that
$G_{\Gamma^{\prime}}$ is the image of the $C$-group
$G_{\widetilde{\Gamma}^{\prime}}$ under the $C$-homomorphism
$i_*:G_{\widetilde{\Gamma}^{\prime}}\to G_{\Gamma}$ given by
embedding $i:\widetilde{\Gamma}^{\prime}\hookrightarrow \Gamma$,
where $\widetilde{\Gamma}^{\prime}$ is the $C$-subgraph of $\Gamma$
obtained from $\Gamma^{\prime}$ after deleting all vertices
$v\not\in \Gamma^{\prime}$ and all edges labeled by the vertices
$v\not\in \Gamma^{\prime}$.

The $C$-group $G_{\Gamma}$ acts on $\Gamma$. Therefore
 the homomorphism $i_*$ defines an action of $G_{\widetilde{\Gamma}^{\prime}}$ on $\Gamma$ leaving fixed
each connected component of $\Gamma$.  It is easy to see that if
$\Gamma^{\prime}$ is ample then $G_{\widetilde{\Gamma}^{\prime}}$
acts transitively on the set of vertices of each connected component
of $\Gamma$.

\begin{lem} \label{ample} Let a union $\Gamma^{\prime}$ of some connected components of a
finite $C$-graph $\Gamma$ is ample and let $v_o$, $v_e$ be two
vertices belonging to a connected component, say $\Gamma_1$, of
$\Gamma$. Then there is a path $l(v_o;v_1, \dots v_n)$ connecting
$v_o$ and $v_e$ such that $v_1,\dots , v_n$ are vertices of
$\Gamma^{\prime}$.
\end{lem}
\proof Let $V_1$ be the set of vertices of $\Gamma_1$. For $v\in
V_1$ denote by $V_1(v)$ the set of all vertices $v^{\prime}\in V_1$
such that for $v^{\prime}$ there is a path $l(v;v_1, \dots v_n)$
connecting $v$ and $v^{\prime}$ such that $v_1,\dots , v_n$ are
vertices of $\Gamma^{\prime}$. We say that $v\leq v^{\prime}$ if
$v^{\prime}\in V_1(v)$. It is easy to see that  $V_1(v)\subset
V_1(v^{\prime})$ if $v^{\prime}\leq v$. Therefore the set of subsets
$V_1(v)\subset V_1$, $v\in V_1$, is partially ordered under
inclusions and hence there is a maximal one, say $V_1(\tilde v)$,
since $\Gamma$ is a finite graph.

Next, it is easy to see that the $C$-group
$G_{\widetilde{\Gamma}^{\prime}}$ acts on the set of subsets
$V_1(v)\subset V_1$, $v\in V_1$, by the rule:
$y_{v_1}(V_1(v))=V_1(y_{v_1}(v))$ for a $C$-generator $y_{v_1}\in
G_{\widetilde{\Gamma}^{\prime}}$. Therefore for a maximal subset
$V_1(\tilde v)$ we have $V_1(\tilde v)=V_1(v)$ for all $v\in
V_1(\tilde v)$, since $G_{\widetilde{\Gamma}^{\prime}}$ acts
transitively on $V_1$. Hence, $V_1$ can be represented as the
disjoint union $V_1(\tilde v_1)\sqcup \dots \sqcup V_1(\tilde v_m)$
of maximal subsets $V_1(\tilde v_i)$. Finally,
since $\Gamma^{\prime}$ is ample, we obtain that $m=1$. \qed \\

Let a union $\Gamma^{\prime}$ of some connected components of a
finite $C$-graph $\Gamma$ be ample and let $\Gamma_1$ be a connected
component of $\Gamma$. By definition, the {\it distance}
$d_{\Gamma^{\prime}}(v_o,v_e)$ {\it between two vertices} $v_o$,
$v_e$ of $\Gamma_1$ {\it with respect to} $\Gamma^{\prime}$ is the
smallest $n$ such that there is a path $l(v_o;v_1, \dots v_n)$
connecting $v_o$ and $v_e$ such that $v_1,\dots , v_n$ are vertices
of $\Gamma^{\prime}$. The number
$d_{\Gamma^{\prime}}(\Gamma_1)=\max_{v_o,v_e\in
V_1}d_{\Gamma^{\prime}}(v_o,v_e)$ is called the {\it diameter} of
$\Gamma_1$ {\it with respect to} $\Gamma^{\prime}$.

\begin{prop}\label{rep2} Let $G_{\Gamma}=(G,Y)$ be a finite
$C$-group, $\Gamma^{\prime}=\Gamma_1\sqcup \dots  \sqcup \Gamma_k$ a
ample subgraph of $\Gamma$, and $\Gamma\setminus
\Gamma^{\prime}=\Gamma_{k+1}\sqcup\dots\sqcup \Gamma_m$, where
$\Gamma_i$, $i=1,\dots, m$, are the connected components of
$\Gamma$. Denote by
$d=d_{\Gamma^{\prime}}=\max(d_{\Gamma^{\prime}}(\Gamma_{1}),\dots ,
d_{\Gamma^{\prime}}(\Gamma_m))$. Let an element $s\in S(G,Y)^G$ is
such that $\tau_i(s)\geq 2n_ip_id+1$ for all $i\leq k$, where
$\tau_i(s)$is the $i$th coordinate of $\tau(s)$. Then the element
$s$ can be written in the form: $s=(x_{v_{{k+1},1}}^{a_{k+1}}\cdot
\, .\, .\, .\, \cdot x_{v_{m,1}}^{a_m})\cdot
s_{\Gamma^{\prime}}^d\cdot s_1$, where $a_i=\tau_i(s)$ for
$i=k+1,\dots, m$ and $s_1\in S(G,Y)$ is such that $\tau_i(s_1)=0$
for $i=k+1,\dots, m$.
\end{prop}
\proof Let us write the element $s$ in the form: $s=s^{\prime}\cdot
s^{\prime\prime}$ where $s^{\prime}$ and $s^{\prime\prime}$ are such
that $\tau_i(s^{\prime})=0$ for $i\geq k+1$ and
$\tau_i(s^{\prime\prime})=0$ for $i\leq k$. Then, by Corollary
\ref{can}, the element $s$ can be written in the form: $s=\widetilde
s^{\prime}\cdot s^d_{\Gamma^{\prime}}\cdot s^{\prime\prime}$.

Let for some $j$ a letter $x_{v_{m,j}}$ enter in $s^{\prime\prime}$.
Connect the vertex $v_{m,1}$ with $v_{m,j}$ by a path
$l(v_{m,1};v_1,\dots, v_r)$ of length $r\leq d$ (remind that
$\Gamma^{\prime}$ is ample, therefore by Lemma \ref{ample},
$v_{m,1}$ and $v_{m,j}$ can be connected by such a path), where
$v_1, \dots, v_r$ are some vertices of $\Gamma^{\prime}$, and write
$s^{\prime\prime}$ in the form $s^{\prime\prime}=x_{v_{m,j}}\cdot
\widetilde s^{\prime\prime}$. By Lemma \ref{divis},
$s^d_{\Gamma^{\prime}}=\widetilde s\cdot (x_{v_1}\cdot\, .\, .\, .
\cdot x_{v_r})$ for some $\widetilde s$. We have
$\alpha(x_{v_1}\cdot\, .\, .\, . \cdot x_{v_r})=y_{v_1}\dots
y_{v_r}$ and by definition of $C$-graphs of $C$-groups, we have
$(y_{v_1}\dots y_{v_r})y_{v_{m,j}}(y_{v_1}\dots
y_{v_r})^{-1}=y_{v_{m,1}}$. Therefore,
$$s^d_{\Gamma^{\prime}}\cdot x_{v_{m,j}}=\widetilde s\cdot (x_{v_1}\cdot\, .\, .\, .
\cdot x_{v_r})\cdot x_{v_{m,j}}=\widetilde s\cdot x_{v_{m,1}}\cdot
(x_{v_1}\cdot\, .\, .\, . \cdot x_{v_r})$$ and hence (after moving
$x_{v_{m,1}}$ to the right)
$$s=s_1^{\prime}\cdot s^d_{\Gamma^{\prime}}\cdot  s^{\prime\prime}=
\widetilde s^{\prime}\cdot \widetilde s\cdot  x_{v_{m,1}}\cdot
(x_{v_1}\cdot\, .\, .\, . \cdot x_{v_r})\cdot \widetilde
s^{\prime\prime}= s_1^{\prime}\cdot s_1^{\prime\prime}\cdot
x_{v_{m,1}},$$ where $s_1^{\prime}$ is such that
$\tau_i(s_1^{\prime})=0$ for $i\geq k+1$ and
$\tau_i(s_1^{\prime})\geq 2n_ip_id+1$ for all $i\leq k$, and
$s_1^{\prime\prime}$ is such that  $\tau_i(s_1^{\prime\prime})=0$
for $i\leq k$.

If again for some $j$ a letter $x_{v_{m,j}}$ enters in
$s_1^{\prime\prime}$, applying Corollary \ref{can},  Claim
\ref{ample}, and Lemma \ref{divis}, we can repeat the transformation
described above, and so on, and we obtain a factorization $s=
s_{a_m}^{\prime}\cdot s_{a_m}^{\prime\prime}\cdot
x^{a_m}_{v_{m,1}}$. After that we can repeat the transformation
described above and we obtain a factorization $s=
s_{a_m+a_{m-1}}^{\prime}\cdot s_{a_m+a_{m-1}}^{\prime\prime}\cdot
x^{a_{m-1}}_{v_{m-1,1}}\cdot x^{a_{m}}_{v_{m,1}}$, and so on. After
the last step of these transformations, we move the obtained product
$(x_{v_{k+1,1}}^{a_{k+1}}\cdot \, .\, .\, .\, \cdot
x_{v_{m,1}}^{a_m})$ to the left and apply Corollary \ref{can} to
complete the proof. \qed
\subsection{$C$-graphs of equipped groups} To each equipped group
$(G,O)$, one can associate a $C$-graph. By definition, the $C$-{\it
graph $\Gamma=\Gamma_{(G,O)}$ of an equipped group} $(G,O)$ is the
$C$-graph whose set of vertices is in one to one correspondence with
the set $O$. Two vertices $v_{g_1}$ and $v_{g_2}$, $g_1, g_2\in O$
are connected by the labeled edge $e_{v_{g_1},v_{g_2},v_g}$ if and
only if $g^{-1}g_1g=g_2$ for some $g\in O$.

\begin{ex} {\rm To describe the $C$-graph $\Gamma_{(\Sigma_n,T_n)}$ of
the equipped symmetric group $(\Sigma_n,T_n)$, where $T_n$ is the
set of transpositions, consider a $(n-1)$-simplex $\Delta_{n-1}$.
Let $V_1,\dots V_n$ be the vertices of $\Delta_{n-1}$ and $E_{i,j}$
its edges connecting the vertices $V_i$ and $V_j$. The vertices
$v_{i,j}$ of $\Gamma_{(\Sigma_n,T_n)}$ are the middles of the edges
$E_{i,j}$. If $E_{i,j}$ and $E_{k,l}$ are skewed edges, then the
edge of $\Gamma_{(\Sigma_n,T_n)}$ with tail $v_{i,j}$ and label
$v_{k,l}$ is the loop. The head of the edge with tail $v_{i,j}$ and
label $v_{j,k}$ is $v_{i,k}$.}
\end{ex}

To each equipped group $(G,O)$, we can associate a $C$-group
$G_{\Gamma_{(G,O)}}$. Denote the $C$-generators $y_{v_{g_i}}$,
$g_i\in O$, of the $C$-group $G_{\Gamma_{(G,O)}}$ by $y_{g_i}$. We
have the natural homomorphism of equipped groups
$\beta=\beta_{(G,O)}:G_{\Gamma_{(G,O)}}\to (G,O)$ given by
$\beta(y_{g_i})=g_i$ for all $g_i\in O$. Obviously, $\beta_{|Y}:
Y\to O$ is one to one correspondence.
\begin{claim} \label{beta} {\rm ($i$)} $\ker \beta$ is a subgroup of the center of $G_{\Gamma_{(G,O)}}$.
\newline {\rm ($ii$)} If $(G,O)$ is an equipped group such that $G$
is generated by the elements of $O$, then $\beta$ and $\beta_{\mid
[G_{\Gamma_{(G,O)}},G_{\Gamma_{(G,O)}}]}:
[G_{\Gamma_{(G,O)}},G_{\Gamma_{(G,O)}}]\to [G,G]$ are epimorphisms.
\end{claim}
\proof Obvious. \qed

\subsection{Equivalence of equipped groups}\label{amb}  Let $(G_1,O_1)$ and $(G_2,O_2)$ be two equipped groups
such that $G_1$ and $G_2$ are generated, resp., by the elements of
$O_1$ and $O_2$. We say that $(G_1,O_1)$ and $(G_2,O_2)$ are {\it
equivalent} if the $C$-graphs $\Gamma_{(G_1,O_1)}$ and
$\Gamma_{(G_2,O_2)}$ are isomorphic as $C$-graphs.

\begin{claim} \label{iso3} Let $(G_1,O_1)$ and $(G_2,O_2)$ be two equivalent equipped
groups. Then  the $C$-groups $G_{\Gamma_{(G_1,O_1)}}$ and
$G_{\Gamma_{(G_2,O_2)}}$ are $C$-isomorphic.
\end{claim}
\proof It follows from Claim \ref{iso4}. \qed \\

It follows from Claims \ref{beta} and \ref{iso4} that for each class
of equivalent equipped groups corresponding to a $C$-graph $\Gamma$,
there is a maximal one, namely, the $C$-group $G_{\Gamma}$, such
that for any equipped group $(G,O)$ belonging to this class there is
an epimorphism of equipped groups, namely
$\beta_{(G,O)}:G_{\Gamma_{(G,O)}}\to (G,O)$, which is defined
uniquely by an isomorphism  $\Gamma\simeq \Gamma_{(G,O)}$ and by the
following condition: for each $C$-generator $y$ of $G_{\Gamma}$ the
image $\beta_{(G,O)}(y)=g\in O$ if $y$ and $g$ correspond to the
same vertex $v$ of $\Gamma\simeq \Gamma_{(G,O)}$. By Claim
\ref{beta}, $\ker \beta_{(G,O)}$ is a subgroup of the center
$Z(G_{\Gamma})$ of $G_{\Gamma}$. The inverse statement is also true,
namely, an equipped group $(G,O)$, obtained as the quotient group
$G_{\Gamma}/H$ of a $C$-group $G_{\Gamma}$, is equivalent to
$G_{\Gamma}$ if $H$ is a subgroup of $Z(G_{\Gamma})$ and it contains
neither $C$-generators of $G_{\Gamma}$ nor quotients $y_iy_j^{-1}$
of $C$-generators $y_i$ and $y_j$, $y_i\neq y_j$.

The order $a_{(G,O)}=|H\cap [G_{\Gamma},G_{\Gamma}]|$ of the group
$H\cap [G_{\Gamma},G_{\Gamma}]$ is called the {\it ambiguity index}
of the equipped group $(G=G_{\Gamma}/H,O)$ equivalent to
$G_{\Gamma}$.

\begin{prop} Let an equipped group $(G,O)$ is equivalent to
a $C$-group $G_{\Gamma}$. If $G$ is a perfect group, then
$G_{\Gamma}$ is isomorphic to the direct product
$[G_{\Gamma},G_{\Gamma}] \times
(G_{\Gamma}/[G_{\Gamma},G_{\Gamma}])$.
\end{prop}
\proof For each connected component $\Gamma_i$ of $\Gamma$, let us
choose a vertex $v_i\in \Gamma_i$.

The restriction of $\beta_{(G,O)}$ to $[G_{\Gamma},G_{\Gamma}]$ is
an epimorphism onto $G$, since $\beta_{(G,O)}$ is epimorphism and
$G$ is a perfect group. Therefore for each $C$-generator $y_{v_i}$
of $G_{\Gamma}$, there is an element $g_i\in
[G_{\Gamma},G_{\Gamma}]$ such that
$\beta_{(G,O)}(g_i)=\beta_{(G,O)}(y_{v_i})$. We have
$y_{v_i}g_i^{-1}\in \ker \beta_{(G,O)}\subset Z(G_{\Gamma})$ and
$ab(y_{v_i}g_i^{-1})=(0,\dots ,0,1,0,\dots, 0)$, where $1$ stands on
the $i$th place. Therefore the elements $y_{v_i}g_i^{-1}$ generate
in $Z(G_{\Gamma})$ a free abelian group $H$ such that $ab_{\mid H}:
H\to G_{\Gamma}/[G_{\Gamma},G_{\Gamma}]$ is an isomorphism and
Proposition follows from the short exact sequence
\begin{equation}
\label{exac} 1\to [G_{\Gamma},G_{\Gamma}]\to G_{\Gamma}\to
G_{\Gamma}/[G_{\Gamma},G_{\Gamma}]\to 1. \end{equation}

\begin{prop} \label{pr} Let an equipped group $(G,O)$ is equivalent to
a $C$-group $G_{\Gamma}$. If $O$ consists of a single conjugasy
class, then the group $G_{\Gamma}$ is isomorphic to the semidirect
product $[G_{\Gamma},G_{\Gamma}] \rtimes \mathbb Z. $
\end{prop}
\proof It follows from exact sequence (\ref{exac}), since
$G_{\Gamma}/[G_{\Gamma},G_{\Gamma}]\simeq \mathbb Z$ if $O$ consists
of a single conjugacy class. \qed \\

Note that Proposition \ref{pr} is not true if $O$ consists of more
than one conjugacy class. For example, Proposition \ref{pr} is not
true if $G_{\Gamma}$ is a free group $\mathbb F^n$, $n>1$, whose set
of $C$-generators is the union of conjugacy classes of a set of free
generators of $\mathbb F^n$.

\begin{lem} \label{alg1} Let a subgroup $H$ of the center $Z(G_{\Gamma})$
of a finite $C$-group $G_{\Gamma}$ be generated by the elements
$y_{v_i}^{k_ip_i}\in Z(G_{\Gamma})$, $i=1,\dots,m$, where
$k_ip_i\geq 2$, $v_i$ is a vertex of the $i$th connected component
of the $C$-graph $\Gamma$ and $p_i$ is its period. Then
$(G,O)=G_{\Gamma}/H$ is an equipped group equivalent to $G_{\Gamma}$
and $\beta_{(G,O)}: [G_{\Gamma_{(G,O)}},G_{\Gamma_{(G,O)}}]\to
[G,G]$ is an isomorphism. In particular,  $a_{(G,O)}=1$.
\end{lem}
\proof Obvious. \qed
\subsection{The type homomorphism} \label{type}
Let $(G,O)$ be an equipped group, $\Gamma=\Gamma_{(G,O)}$ its
$C$-graph, and $G_{\Gamma}$ 
the $C$-group equivalent to $(G,O)$. The homomorphism
$\beta=\beta_{(G,O)}: G_{\Gamma}\to (G,O)$ defines a homomorphism of
semigroups $\beta_*:S(G_{\Gamma},Y)\to S(G,O)$ given by
$\beta_*(x_{y_g})=x_{g}$ for $g\in O$.
\begin{claim}\label{is} The homomorphism $\beta_*$ is an isomorphism.
\end{claim}
\proof Obvious.

In the sequel, according to Claim \ref{is}, we will identify the
semigroups $S(G,O)$ and $S(G_{\Gamma_{(G,O)}},Y)$.

The homomorphism of semigroups $\tau :
S(G,O)=S(G_{\Gamma_{(G,O)}},Y)\to \mathbb Z_{\geq 0}^m\subset
\mathbb
Z^m=G_{\Gamma_{(G,O)}}/[G_{\Gamma_{(G,O)}},G_{\Gamma_{(G,O)}}]$ is
called the {\it type homomorphism} and the image $\tau(s)$ of $s\in
S(G,O)$ is called the {\it type} of $s$. If $O$ consists of a single
conjugacy class, then the homomorphism $\tau$ can (and will) be
identified with the homomorphism $ln:S(G,O)\to \mathbb Z_{\geq 0}$.
In general case, if $\tau(s)=(\tau_1(s),\dots , \tau_m(s))$, then
$ln(s)=\sum_{i=1}^m\tau_i(s)$.

An element $g\in G_{\Gamma}$ is called {\it positive} if there is
$s\in S(G_{\Gamma},Y)$ such that $\alpha_{G_{\Gamma}}(s)= g$.
\begin{lem} \label{simple2} Any element $g$ of the $C$-group $G_{\Gamma}$
can be represented in the form:
\begin{equation} \label{diferen} g=g_1
g_2^{-1},\end{equation} where $g_1$ are $g_2$ are positive elements.
In particular, $g\in [G_{\Gamma},G_{\Gamma}]$ if and only if
$ab(g_1)=ab(g_2)$ in representation {\rm (\ref{diferen})} of $g$ as
a quotient of two positive elements $g_1$ and $g_2$.

If $G_{\Gamma}$ is a finite $C$-group, then for each $g\in
G_{\Gamma}$ there is a presentation {\rm (\ref{diferen})} of $g$ as
a quotient of two positive elements $g_1$ and $g_2$ such that $g_2=
\alpha_{G_{\Gamma}}(s_{\Gamma}^n)$ for some $n$ and $g_1=
\alpha_{G_{\Gamma}}(s_{1})$ for some $s_1\in
S(G_{\Gamma},Y)^{G_{\Gamma}}$, where $s_{\Gamma}$ is the canonical
element of $S(G_{\Gamma},Y)$.
\end{lem}
\proof The first part of Lemma \ref{simple2} is obvious. The second
part follows from Lemma \ref{divis}. \qed

Let $H$ be a subgroup of the center of a finite $C$-group
$G_{\Gamma}$. Denote by $z_H$ the minimal exponent $n$ such that
each element $g\in H\cap [G_{\Gamma},G_{\Gamma}]$ has   presentation
(\ref{diferen}) in which $g_2=\alpha_{G_{\Gamma}}(s_{\Gamma}^{n})$.
\begin{prop}\label{rep5} Let $G_{\Gamma}=(G,Y)$ be a finite
$C$-group, $\Gamma^{\prime}=\Gamma_1\sqcup \dots  \sqcup
\Gamma_k\subset \Gamma$ be an ample subgraph of $\Gamma$, and
$\widetilde{\Gamma^{\prime}}$ the $C$-graph defined by
$\Gamma^{\prime}$. Then ${i_*}_{\mid
[G_{\widetilde{\Gamma^{\prime}}},G_{\widetilde{\Gamma^{\prime}}}]}:
[G_{\widetilde{\Gamma^{\prime}}},G_{\widetilde{\Gamma^{\prime}}}]
\to [G_{\Gamma},G_{\Gamma}]$ is an epimorphism.
\end{prop}
\proof It follows from  Proposition \ref{rep2} and Lemma
\ref{simple2}. \qed
\begin{cor} Let $G_{\Gamma}$ and $G_{\widetilde{\Gamma}^{\prime}}$ be finite $C$-groups
equivalent respectively to equipped groups $(G,O)$ and
$(G,O^{\prime})$, where $O^{\prime}\subset O$ are equipments of $G$
such that the elements of $O^{\prime}$ generate the group $G$. Then
$a_{(G,O)}\leq a_{(G,O^{\prime})}$, in particular, if
$a_{(G,O^{\prime})}=1$, then $a_{(G,O)}=1$.
\end{cor}

\section{Stability of factorization semigroups}
\subsection{Equivalence of elements}
Elements $s_1$ and $s_2$ of a semigroup $S$  are said to be {\it
$r$-equivalent} (resp., {\it $l$-equivalent}) if there is an element
$s_3\in S$ such that $s_1\cdot s_3=s_2\cdot s_3$ (resp., $s_3\cdot
s_1=s_3\cdot s_2$), and they are {\it equivalent} if there are two
elements $s_3$, $s_4\in S$ such that $s_3\cdot s_1\cdot s_4=
s_3\cdot s_2\cdot s_4$. Notation $s_1\sim s_2$ (resp., $s_1\sim_r
s_2$ and $s_1\sim_l s_2$) means that elements $s_1$ and $s_2$ are
equivalent (resp., $r$-equivalent and $l$-equivalent). It is easy to
see that if $s_1\sim s_2$, where $s_1, s_2\in S(G,O)$, then
$\tau(s_1)=\tau(s_2)$.

\begin{lem} \label{eqival} Let $s_1, s_2\in S(G,O)$ be two elements of a factorization
semigroup over a group $G$. Then the following statements are
equivalent:
\begin{itemize}
\item[$(i)$] $s_1\sim_r s_2$;
\item[$(ii)$] $s_1\sim_l s_2$;
\item[$(iii)$] $s_1\sim s_2$.
\end{itemize}
\end{lem}
\proof We prove only the implication $(i)\Rightarrow (ii)$, since
the proof of all other implications are similar.

We have $s_1\cdot s_3=\rho(\alpha_G(s_1))(s_3)\cdot s_1$ and
$s_2\cdot s_3=\rho(\alpha_G(s_2))(s_3)\cdot s_2$. Since $\alpha_G$
is a homomorphism, it is easy to see that $\alpha_G
(s_1)=\alpha_G(s_2)$ if $s_1$ and $s_2$ are $r$-equivalent (resp.,
$l$-equivalent or equivalent). Therefore, if $s_1\sim_r s_2$, that
is, $s_1\cdot s_3=s_2\cdot s_3$ for some $s_3\in S(G,O)$, then
$\rho(\alpha_G(s_1))(s_3)\cdot s_1=\rho(\alpha_G(s_2))(s_3)\cdot
s_2$, that is, $s_1\sim_l s_2$.
\qed \\

\begin{lem} \label{eqival2} Let $S(G,O)$ be a factorization
semigroup over a group $G$. Then the relation $s_1\sim s_2$ is an
equivalence relation.
\end{lem}
\proof Let we have $s_1\sim s_2$ and $s_2\sim s_3$. Then, by Lemma
\ref{eqival}, there are elements $s_4$ and $s_5$ such that $s_4\cdot
s_1=s_4\cdot s_2$ and $s_2\cdot s_5=s_3\cdot s_5$. Therefore
$s_4\cdot s_1\cdot s_5=s_4\cdot s_2\cdot s_5=s_4\cdot s_3\cdot s_5$,
that is, $s_1\sim s_3$.
\qed \\
\begin{lem} \label{eqival3} Let $s_1, s_2,s_3,s_4\in S(G,O)$ be four elements of a factorization
semigroup over a group $G$. Then $s_1\cdot s_2\sim s_3\cdot s_4$ if
$s_1\sim s_3$ and $s_2\sim s_4$. If $s_1\cdot s_3\sim s_2\cdot s_4$
and $s_1\sim s_2$, then $s_3\sim s_4$.
\end{lem}
\proof Similar to the proof of Lemma \ref{eqival2}. \qed \\

\begin{thm} \label{main1} Let $G$ be a $C$-group and $Y$ the set of its $C$-generators.
Two elements $s_1$ and $s_2\in S(G,Y)$ are equivalent if and only if
$\alpha_G(s_1)=\alpha_G(s_2)$, where $\alpha_G$ is the product
homomorphism.
\end{thm}
\proof It is obvious that if $s_1\sim s_2$, then
$\alpha_G(s_1)=\alpha_G(s_2)$.

Denote by $x_i=x_{y_i}$ the generator of $S(G,Y)$ corresponding to
the $C$-generator $y_i\in Y$. Let $s_1=x_{i_1}\cdot \, .\, .\, .\,
\cdot x_{i_n}$ and $s_2=x_{j_1}\cdot \, .\, .\, .\, \cdot x_{j_k}$
be such that $\alpha_G(s_1)=y_{i_1}\dots y_{i_n}=y_{j_1}\dots
y_{j_m}=\alpha_G(s_2)$. Then it is easy to see that
$\tau(s_1)=\tau(s_2)$ and, in particular, $n=k$. In addition, the
word $w=y_{i_1}\dots y_{i_n}y_{j_n}^{-1}\dots y_{j_1}^{-1}$ in
letters of $Y$ represents the unity of $G$.

To prove $s_1\sim s_2$ we will use some admissible transformations
of  van Kampen diagrams defined over the $C$-presentation of $G$. To
define them, recall that by van Kampen Lemma (see, for example,
\cite{L-S} or \cite{O}), for the word $w$ there is a van Kampen
diagram, that is, a planar finite cell complex $D\subset \mathbb
R^2$ with the following additional data and satisfying the following
additional properties:

1. The complex $D$ is connected and simply connected.

2. Each edge (one-cell) of $D$ is directed and labeled by a letter
$y\in Y$.

3. Some vertex (zero-cell) which belongs to the topological boundary
$\partial D$ of $D$  is specified as a base-vertex $A$ called the
{\it origin} of the diagram.

4. Each region (two-cell) of $D$ is a quadrangle corresponding to a
$C$-relation $y_iy_j=y_ky_i$ of the $C$-presentation of $G$ (see
Fig. $1$; the vertex $v_1$ will be called the {\it bottom} of $Q$
and the vertex $v_2$ will be called the {\it top} of $Q$).

\begin{picture}(0,90)
\put(224,65){$v_2$}\put(220,65){\circle*{3}}
\put(220,25){\circle*{3}}\put(180,25){\vector(1,0){40}}
\put(200,17){$y_i$}
\put(220,25){\vector(0,1){40}}\put(224,45){$y_j$}
\put(180,25){\vector(0,1){40}}\put(168,45){$y_k$}
\put(173,18){$v_1$}\put(180,65){\circle*{3}}
\put(180,25){\circle*{3}}\put(180,65){\vector(1,0){40}}
\put(200,68){$y_i$} \put(195,0){$\text{Fig.} 1$}
\end{picture} \\

5. The boundary cycle $\partial D$, that is, an edge-path
corresponding to going around once in the clockwise direction along
the boundary of the unbounded complementary region of $D$, starting
and ending at the origin $A$, has the label $w=y_{i_1}\dots
y_{i_n}y_{j_n}^{-1}\dots y_{j_1}^{-1}$.

Since in our case the word $w$ splits into two subwords: the first
one consists of letters with positive exponents and the other one
consists of letters with negative exponents, the van Kampen diagram
$D$ is of the following form: it is a chain of discs connected by
simple directed paths (see Fig. 2).

\begin{picture}(0,90)
\put(80,40){\vector(1,0){20}}
\put(100,40){\vector(1,0){20}}\put(120,40){\circle*{3}}\put(160,40){\circle*{3}}
\put(300,40){\circle*{3}}\put(260,40){\circle*{3}}
\put(140,40){\circle{40}}\put(140,60){\vector(1,0){0.1}}\put(140,20){\vector(1,0){0.1}}
\put(160,40){\vector(1,0){20}}
\put(180,40){\vector(1,0){20}}\put(203,40){$\dots$}
\put(220,40){\vector(1,0){20}} \put(240,40){\vector(1,0){20}}
\put(280,40){\circle{40}}\put(280,60){\vector(1,0){0.1}}\put(280,20){\vector(1,0){0.1}}
\put(300,40){\vector(1,0){20}}
\put(320,40){\vector(1,0){20}}\put(342,36){$B$}
\put(70,36){$A$}\put(80,40){\circle*{3}}\put(340,40){\circle*{3}}
\put(195,0){$\text{Fig.} 2$}
\end{picture} \\ \\
Denote the end vertex of this chain by $B$ and call it the {\it end}
of $D$. Diagrams satisfying conditions 1 -- 5 will be called {\it
admissible.}

It follows from Lemmas \ref{eqival2} and \ref{eqival3} that it
suffices to consider only the case when $D$ consists of a single
disc (see Fig. $3$). \\ \\

\begin{picture}(0,100)
\put(185,15){$y_{i_1}$}\put(185,82){$y_{i_n}$}
\put(203,15){$y_{j_1}$}\put(208,82){$y_{j_n}$}
\put(120,10){\vector(0,1){20}}\put(120,30){\vector(0,1){20}}\put(120,70){\vector(0,1){20}}
\put(144,10){$\dots$}\put(144,90){$\dots$}
\put(120,53){$\cdot$}\put(120,57){$\cdot$}\put(120,62){$\cdot$}
\put(120,90){\vector(1,0){20}}\put(140,10){\vector(-1,0){20}}
\put(120,90){\circle*{3}}\put(120,10){\circle*{3}}
\put(160,90){\vector(1,0){20}}\put(180,90){\vector(1,0){20}}\put(220,90){\vector(-1,0){20}}
\put(240,90){\vector(-1,0){20}}\put(280,90){\vector(-1,0){20}}
\put(245,10){$\dots$}\put(245,90){$\dots$}
\put(280,90){\circle*{3}}\put(280,30){\vector(0,1){20}}
\put(280,53){$\cdot$}\put(280,57){$\cdot$}\put(280,62){$\cdot$}\put(280,70){\vector(0,1){20}}
\put(200,10){\vector(1,0){20}}\put(200,10){\vector(-1,0){20}}\put(180,10){\vector(-1,0){20}}
\put(260,10){\vector(1,0){20}}\put(280,10){\vector(0,1){20}}\put(280,10){\circle*{3}}
\put(220,10){\vector(1,0){20}}\put(195,95){$B$}
\put(195,-5){$A$}\put(200,10){\circle*{3}}\put(200,90){\circle*{3}}
\put(195,-30){$\text{Fig.} 3$}
\end{picture} \\ \\ \\

A vertex $v$ in a van Kampen diagram $D$ is called {\it locally
maximal} (resp., {\it minimal}), if there is not an edge of $D$ for
which $v$ is the tail (resp., head). A path $l$ along edges $e_1,
\dots , e_k$ of $D$ is called {\it increasing} if the tail of each
edge $e_{i+1}$ is the head of the edge $e_i$ for $i=1,\dots, k$. The
label of increasing path $l$ is a positive word $y_{l_1}\dots
y_{l_k}$.

The origin $A$ and the end $B$ divides the boundary $\partial D$
into two parts. The increasing path $\partial_l D$ (resp.,
$\partial_r D$) along the boundary $\partial D$ connecting $A$ and
$B$ and having the label $y_{i_1}\dots y_{i_n}$ (resp.,
$y_{j_1}\dots y_{j_n}$) will be called the {\it left} (resp., {\it
right}) {\it side} of $\partial D$.

If the origin $A$ is not a locally minimal vertex of $D$, then there
are a locally minimal vertex $C$ and a simple increasing path $l$
connecting the vertices $C$ and $A$. (Note that the vertex $C$ can
not be a vertex belonging to the boundary $\partial D$, since a
positive word in $C$-generators can not represent the unity of a
$C$-group.) Then we can cut the disc $D$ along the path $l$ and, as
a result, we obtain a new disc diagram $D^{\prime}$ (see Fig. $4$)
in which $C$ is the origin and the label of $\partial D^{\prime}$ is
$$y_{l_1}\dots y_{l_k}y_{i_1}\dots y_{i_n}y^{-1}_{j_n}\dots y^{-1}_{j_1} y_{l_k}^{-1}\dots y_{l_1}^{-1},$$
where $y_{l_1}\dots y_{l_k}$ is the label of $l$.
 \\ \\ \\

 \begin{picture}(0,100)
\put(65,15){$y_{i_1}$} \put(59,92){$y_{i_n}$}
\put(83,15){$y_{j_1}$}\put(88,92){$y_{j_n}$}
\put(0,10){\vector(0,1){20}}\put(0,30){\vector(0,1){20}}\put(0,80){\vector(0,1){20}}
\put(24,10){$\dots$}\put(24,100){$\dots$}
\put(0,55){$\cdot$}\put(0,59){$\cdot$}\put(0,65){$\cdot$}
\put(0,100){\vector(1,0){20}}\put(20,10){\vector(-1,0){20}}
\put(0,100){\circle*{3}}\put(0,10){\circle*{3}}
\put(40,100){\vector(1,0){20}}\put(60,100){\vector(1,0){20}}\put(100,100){\vector(-1,0){20}}
\put(120,100){\vector(-1,0){20}}\put(160,100){\vector(-1,0){20}}
\put(125,10){$\dots$}\put(125,100){$\dots$}
\put(160,100){\circle*{3}}\put(160,30){\vector(0,1){20}}
\put(160,55){$\cdot$}\put(160,59){$\cdot$}\put(160,65){$\cdot$}\put(160,80){\vector(0,1){20}}
\put(80,10){\vector(1,0){20}}\put(80,10){\vector(-1,0){20}}\put(60,10){\vector(-1,0){20}}
\put(140,10){\vector(1,0){20}}\put(160,10){\vector(0,1){20}}\put(160,10){\circle*{3}}
\put(100,10){\vector(1,0){20}}\put(75,105){$B$}
\put(75,-5){$A$}\put(80,10){\circle*{3}}\put(80,100){\circle*{3}}
\put(80,75){$C$}\put(80,70){\circle*{3}}\put(80,70){\vector(0,-1){20}}
\put(83,59){$y_{l_1}$}\put(80,35){\vector(0,-1){25}}\put(83,29){$y_{l_k}$}
\put(80,37){$\cdot$}\put(80,41){$\cdot$}\put(80,44){$\cdot$}
\put(180,45){$\rightsquigarrow$} \put(275,15){$y_{l_1}$}
\put(269,92){$y_{i_n}$}
\put(293,15){$y_{l_1}$}\put(298,92){$y_{j_n}$}
\put(217,14){$y_{l_k}$}\put(196,16){$y_{i_1}$}
\put(353,15){$y_{l_k}$}\put(372,16){$y_{j_1}$}
\put(210,10){\vector(0,1){20}}\put(210,30){\vector(0,1){20}}\put(210,80){\vector(0,1){20}}
\put(234,10){$\dots$}\put(234,100){$\dots$}
\put(210,55){$\cdot$}\put(210,59){$\cdot$}\put(210,65){$\cdot$}
\put(210,100){\vector(1,0){20}} \put(230,10){\vector(-1,0){20}}
\put(210,100){\circle*{3}}\put(210,10){\circle*{3}}
\put(250,100){\vector(1,0){20}}\put(270,100){\vector(1,0){20}}\put(310,100){\vector(-1,0){20}}
\put(330,100){\vector(-1,0){20}} \put(370,100){\vector(-1,0){20}}
\put(335,10){$\dots$}\put(335,100){$\dots$}
\put(370,100){\circle*{3}}\put(370,30){\vector(0,1){20}}
\put(370,55){$\cdot$}\put(370,59){$\cdot$}\put(370,65){$\cdot$}\put(370,80){\vector(0,1){20}}
\put(290,10){\vector(1,0){20}}\put(290,10){\vector(-1,0){20}}
\put(350,10){\vector(1,0){20}}\put(370,10){\vector(0,1){20}}\put(370,10){\circle*{3}}
\put(310,10){\vector(1,0){20}}\put(285,105){$B$}
\put(285,-5){$C$}\put(290,10){\circle*{3}}\put(290,100){\circle*{3}}
\put(207,-5){$A^{\prime}$}\put(365,-5){$A^{\prime\prime}$}
\put(270,10){\vector(-1,0){20}} \put(180,-30){$\text{Fig.} 4$}
\end{picture} \\ \\ \\

We call the transformation  $D\rightsquigarrow D^{\prime}$ described
above an  {\it admissible transformation $I$}. Note that after
transformation $I$ the diagram $D^{\prime}$ is admissible and the
origin $C$ of $D^{\prime}$ is a locally minimal vertex.

Similarly, if the end $B$ is not a locally maximal vertex of $D$,
then there are a locally maximal vertex $C$ and a simple increasing
path $l$ connecting  $B$ and $C$. Then we can cut the disc $D$ along
the path $l$ and, as a result, we obtain a new admissible disc
diagram $D^{\prime}$ in which $C$ is the end and the label of
$\partial D^{\prime}$ is $$y_{i_1}\dots y_{i_n}y_{l_1}\dots
y_{l_k}y_{l_k}^{-1}\dots y_{l_1}^{-1}y^{-1}_{j_n}\dots y^{-1}_{j_1}
,$$ where $y_{l_1}\dots y_{l_k}$ is the label of $l$. We call this
transformation $D\rightsquigarrow D^{\prime}$ an {\it admissible
transformation $II$}. Note that after transformation $II$ the end
$C$ of $D^{\prime}$ is a locally maximal vertex.

Let $v\in \partial_l D$ be a vertex of an admissible disc diagram
$D$ such that $v$ is neither the origin $A$ nor the end $B$ of $D$
and there is an edge $e$ of $D$, $e\not\subset \partial D$, for
which $v$ is the tail. Denote by $y_0$ the label of this edge and
let $C$ be its head, $C\not\in \partial D$. Assume for definiteness
that $v\in \partial_l D$ (the case when $v\in \partial_r D$ is
similar). Let  $y_{i_k}$ (resp., $y_{i_{k+1}}$) be the label of the
edge belonging to $\partial D$ for which $v$ is the head (resp., the
tail). Let us cut the diagram $D$ along $e$ and after that paste
sequentially $n-k$ additional quadrangles as it is depicted in Fig.
5, where the quadrangle, glued at the $l$th step to the cut disc
diagram along the edges labeled by $y_0$ and $y_{i_{k+l}}$,
corresponds to the relation
$y_{i_{k+l}}y_0=y_0y_{i_{k+l}}^{\prime}$, $l=k+1,\dots, n-k$, in the
$C$-group $G$.

\begin{picture}(0,250)
\put(330,20){\vector(0,1){30}}\put(330,20){\circle*{3}}
\put(300,20){\vector(1,0){30}}
\put(330,175){$\cdot$}\put(330,170){$\cdot$}
\put(330,165){$\cdot$}\put(330,160){$\cdot$}\put(330,155){$\cdot$}
\put(330,135){$\cdot$}\put(330,130){$\cdot$}
\put(330,125){$\cdot$}\put(330,120){$\cdot$}\put(330,115){$\cdot$}
\put(330,110){$\cdot$} \put(330,90){$\cdot$}\put(330,85){$\cdot$}
\put(330,80){$\cdot$}\put(330,75){$\cdot$}\put(330,70){$\cdot$}
\put(330,65){$\cdot$}
\put(330,195){\vector(0,1){30}}\put(330,225){\circle*{3}}\put(280,225){$\dots$}\put(330,225){\vector(-1,0){30}}
\put(210,225){\circle*{3}}\put(197,222){$B^{\prime}$}
\put(240,225){\circle*{3}}\put(237,228){$B$}
\put(240,195){\vector(-1,0){30}}\put(195,205){$y_{i_{n}}^{\prime}
$}\put(243,205){$y_{i_{n}}$}
\put(210,195){\vector(0,1){30}}\put(240,195){\vector(0,1){30}}
\put(240,225){\vector(-1,0){30}}\put(270,225){\vector(-1,0){30}}
\put(250,229){$y_{j_n}$}
\put(220,228){$y_{0}$}\put(220,197){$y_{0}$}\put(220,167){$y_{0}$}
\put(210,180){$\cdot$}\put(210,183){$\cdot$} \put(210,187){$\cdot$}
\put(240,180){$\cdot$}\put(240,183){$\cdot$} \put(240,187){$\cdot$}
\put(240,175){\vector(-1,0){30}}\put(183,160){$y_{i_{k+2}}^{\prime}$}\put(243,160){$y_{i_{k+2}}$}
\put(210,145){\vector(0,1){30}}\put(240,145){\vector(0,1){30}}
\put(220,148){$y_{0}$}\put(220,118){$y_{0}$}\put(220,100){$y_{0}$}
\put(240,145){\vector(-1,0){30}}\put(183,130){$y_{i_{k+1}}^{\prime}$}\put(243,130){$y_{i_{k+1}}$}
\put(210,115){\vector(0,1){30}}\put(240,115){\vector(0,1){30}}
\put(210,115){\circle*{3}}\put(200,112){$C$}
\put(245,112){$v^{\prime}$}\put(245,95){$v^{\prime\prime}$}
\put(240,115){\circle*{3}}\put(240,115){\vector(-1,0){28}}
\put(240,100){\vector(-2,1){30}}\put(240,100){\circle*{3}}\put(226,79){$y_{i_k}$}
\put(240,70){\vector(0,1){30}}\put(240,57){$\cdot$}\put(240,61){$\cdot$}
\put(240,53){$\cdot$}\put(6,35){$y_{i_1}$}\put(253,25){$y_{j_1}$}
\put(240,20){\vector(0,1){30}}\put(240,20){\vector(1,0){30}}\put(280,20){$\dots$}
\put(240,20){\circle*{3}} \put(237,10){$A$}
\put(165,105){$\rightsquigarrow$}
\put(110,135){$\cdot$}\put(110,130){$\cdot$}
\put(110,125){$\cdot$}\put(110,120){$\cdot$}\put(110,115){$\cdot$}
\put(110,110){$\cdot$} \put(110,90){$\cdot$}\put(110,85){$\cdot$}
\put(110,80){$\cdot$}\put(110,75){$\cdot$}\put(110,70){$\cdot$}
\put(110,65){$\cdot$}
\put(80,20){\vector(1,0){30}}\put(110,20){\vector(0,1){30}}\put(110,20){\circle*{3}}
\put(110,150){\vector(0,1){30}}\put(110,180){\circle*{3}}
\put(110,180){\vector(-1,0){30}}\put(60,180){$\dots$}\put(35,185){$y_{j_{n}}$}
\put(50,180){\vector(-1,0){30}}\put(15,185){$B$}\put(20,180){\circle*{3}}
\put(5,160){$y_{i_{n}}$}\put(20,150){\vector(0,1){30}}
\put(20,134){$\cdot$}\put(20,139){$\cdot$}
\put(20,143){$\cdot$}\put(-3,113){$y_{i_{k+1}}$}
\put(20,100){\vector(0,1){30}} \put(11,97){$v$}
\put(54,95){$C$}\put(50,100){\circle*{3}}\put(20,100){\vector(1,0){30}}
\put(33,92){$y_{0}$}\put(20,100){\circle*{3}}\put(6,79){$y_{i_k}$}
\put(20,70){\vector(0,1){30}}\put(20,57){$\cdot$}\put(20,61){$\cdot$}
\put(20,53){$\cdot$}\put(6,35){$y_{i_1}$}\put(33,25){$y_{j_1}$}
\put(20,20){\vector(0,1){30}}\put(20,20){\vector(1,0){30}}\put(60,20){$\dots$}
\put(20,20){\circle*{3}} \put(17,10){$A$} \put(170,-10){$\text{Fig.}
5$}
\end{picture} \\ \\

After gluing these quadrangles we obtain a new admissible disc
diagram $D^{\prime}$  whose end is the vertex $B^{\prime}$ and the
label of $\partial D^{\prime}$ is
$$y_{i_1}\dots y_{i_k}y_0y_{i_{k+1}}^{\prime}\dots y_{i_n}^{\prime}y_{0}^{-1}y_{j_n}^{-1}\dots y^{-1}_{j_1} .$$
We call the transformation $D\rightsquigarrow D^{\prime}$ described
above an {\it admissible transformation $III$ defined by the edge
$e$}.
\begin{claim} \label{equ3} If the elements $x_{i_1}\cdot\, .\, .\,
.\,\cdot x_{i_k}\cdot x_0\cdot x_{i_{k+1}}^{\prime}\cdot \, .\, .\,
.\, \cdot x_{i_{n}}^{\prime}$ and $x_{j_1}\cdot\, .\, .\, .\,\cdot
x_{j_n}\cdot x_0$, obtained from $x_{i_1}\cdot\, .\, .\, .\,\cdot
x_{i_k}\cdot  x_{i_{k+1}}\cdot \, .\, .\, .\, \cdot x_{i_{n}}$ and
$x_{j_1}\cdot\, .\, .\, .\,\cdot x_{j_n}$ after admissible
transformation $III$ defined by the edge $e$, are equivalent, then
$x_{i_1}\cdot\, .\, .\, .\,\cdot x_{i_k}\cdot  x_{i_{k+1}}\cdot \,
.\, .\, .\, \cdot x_{i_{n}}$ and $x_{j_1}\cdot\, .\, .\, .\,\cdot
x_{j_n}$ are also equivalent.
\end{claim}
\proof Obvious. \qed \\

Let $D$ be an admissible disc diagram such that its origin $A$ is a
locally minimal vertex and its end $B$ is a locally maximal vertex.
Let $C$ be a locally maximal vertex of $D$, $C\neq B$. Then
$C\not\in \partial D$, since there is only one locally maximal
vertex belonging to $\partial D$, namely, the end $B$. We say that
$C$ is {\it visible} if there is a vertex $v\in
\partial D$, $v\neq A$, such that $v$ can be connected with $C$ by
increasing path $l$ along edges of $D\setminus \partial D$ and such
that $l\cap \partial D=v$. Let the path $l$ consist of edges
$e_1,\dots, e_k$ with labels $y_1, \dots, y_k$. Perform the sequence
of admissible transformations $III$ defined by the edges $e_1,
\dots, e_k$. As a result, we obtain a new admissible diagram
$D^{\prime}$ in which the origin $A$ is a locally minimal vertex,
the end $B^{\prime}$ is a locally maximal vertex, and the number of
locally maximal vertices is strictly less that the number of locally
maximal vertices of $D$. Such sequence of admissible transformations
will be called a {\it transformation decreasing the number of
locally maximal vertices.}

Let us return to the proof that $s_1\sim s_2$ if
$\alpha_G(s_1)=\alpha_G(s_2)$, where $s_1=x_{i_1}\cdot \, . \, . \,
. \, \cdot x_{i_n}$ and $s_2=x_{j_1}\cdot \, . \, . \, . \, \cdot
x_{j_n}$ are two elements of the factorization semigroup $S(G,Y)$
over a $C$-group $(G,Y)$. Since $\alpha_G(s_1)=\alpha_G(s_2)$, the
word $w=y_{i_1}\dots y_{i_n}y_{j_1}^{-1}\dots y_{j_n}^{-1}$
represents the unity of the group $G$. Therefore, by van Kampen
Lemma, there is a plane disc diagram $D$ over the $C$-presentation
of the group $G$ whose boundary label is the word $w$. By Lemmas
\ref{eqival2} and \ref{eqival3}, as it was mention above, we can
assume that $D$ is a single disc. Conversely, an admissible disc
diagram with boundary label $w=y_{i_1}\dots y_{i_n}y_{j_1}^{-1}\dots
y_{j_n}^{-1}$ defines two elements $s_1(D)=x_{i_1}\cdot \, . \, . \,
. \, \cdot x_{i_n}$ and $s_2(D)=x_{j_1}\cdot \, . \, . \, . \, \cdot
x_{j_n}$ of $S(G,Y)$ such that $\alpha_G(s_1(D))=\alpha_G(s_2(D))$.
Note that if $D \rightsquigarrow D^{\prime}$ is an admissible
transformation $I$ or $II$, or $III$, then, by Lemmas \ref{eqival2},
\ref{eqival3} and Claim \ref{equ3}, $s_1(D)\sim s_2(D)$ if and only
if $s_1(D^{\prime})\sim s_2(D^{\prime})$. Therefore, without loss of
generality, we can assume that $D$ is an admissible disc diagram
such that
\begin{itemize}
\item[$(i)$] the origin $A$ of $D$
is a locally minimal vertex,
\item[$(ii)$] the end $B$ of $D$ is a
locally maximal vertex,
\item[$(iii)$]  there is the only one visible
locally maximal vertex of $D$, namely, the end $B$.
\end{itemize}

Let us show that if an admissible disc diagram $D$ satisfying
conditions $(i)$ -- $(iii)$, then $s_1(D)\sim s_2(D)$. Indeed, if
$D$ consists of a single quadrangle $Q$, then the origin $A$ is the
bottom of $Q$ and the end $B$ is the top of $Q$. Obviously, in this
case we have $s_1(D)= s_2(D)$ (see Fig. $1$).

Now, let $D$ satisfy conditions $(i)$ -- $(iii)$, have $K$ invisible
locally maximal vertices and consist of $k$ quadrangles. Consider
the edge $e_1\subset
\partial_l D$ whose tail is $A$. Let $v_1$ be the head of $e_1$ and $Q$ a
quadrangle such that $e_1\subset
\partial Q$. Then the bottom of $Q$ is the origin $A$, since $A$ is
a locally minimal vertex. Let $v_2$ be the top of $Q$ and
$e_2\subset \partial Q$ the edge connecting $v_1$ and $v_2$.

There are two possibilities: either $e_2\subset
\partial_l D$ and hence $v_2\in \partial_l D$,  or $e_2\not\subset \partial_l D$. In the first case if
$y_{i_1}y_{i_2}y_{l_2}^{-1}y_{l_1}^{-1}$ is the label of $\partial
Q$, then $s_1(D)= x_{i_1}\cdot x_{i_2}\cdot (x_{i_3} \cdot \, .\,
.\, .\, \cdot x_{i_n})=x_{l_1}\cdot x_{l_2}\cdot (x_{i_3} \cdot \,
.\, .\, .\, \cdot x_{i_n})$, and if we cut $Q$ from $D$, then we
obtain a new admissible disc diagram $D^{\prime\prime}=D\setminus Q$
having only $k-1$ quadrangles and such that
$s_1(D^{\prime\prime})=x_{l_1}\cdot x_{l_2}\cdot (x_{i_3} \cdot \,
.\, .\, .\, \cdot x_{i_n})$ and $s_2(D^{\prime\prime})=x_{j_1}\cdot
\, . \, . \, . \, \cdot x_{j_n}$. Denote also $Q$ by $D^{\prime}$.
Note that it is not obligatory that $D^{\prime\prime}$ satisfies
conditions $(i)$ -- $(iii)$. Indeed, if the edge $e$ of $Q$ with
label $x_{l_1}$ belongs to $\partial_r D$, then we must delete it to
obtain admissible disc diagram. Therefore in this case the boundary
label of $\partial D^{\prime\prime}$ is $y_{l_2}(y_{i_3} \dots
y_{i_n})y_{j_n}^{-1}\dots y_{j_2}^{-1}$ (since $y_{l_1}=y_{j_1}$)
and if the origin of $D^{\prime\prime}$ is not a locally minimal
vertex of $D^{\prime\prime}$, then we must perform a transformation
$I$ which increases neither the number of locally maximal vertices
nor the number of quadrangles. Next, it is possible that an
invisible locally maximal vertex of $D$ becomes visible in
$D^{\prime\prime}$. (It is possible only if $K>0$.) But, in this
case we can perform an admissible transformation decreasing the
number of locally maximal vertices and obtain a new admissible disc
diagram having strictly less than $K$ invisible locally maximal
vertices.

If $v_2\not\in \partial D$, then there is an edge $e_3$ in $D$ whose
tail is $v_2$, since, by assumption, $v_2$ is not a locally maximal
vertex. Let $v_3$ is the head of $e_3$. If $v_3\not\in
\partial D$, then there is an edge $e_4$ whose tail is $v_3$, and so
on. As a result, we can find an increasing path $l=(e_2, \dots e_m)$
connecting the vertex $v_1$ and a vertex $v_m$ belonging to
$\partial D$. There are two possibilities: either $v_m\in
\partial_l D$ (see Fig. $6$, Case $l$) or $v_m\in
\partial_r D$ (see Fig. $6$, Case $r$).
In both cases the path $l$ divides $D$ into two disc diagrams
$D^{\prime}$ and $D^{\prime\prime}$ each of which consists not more
than $k-1$ quadrangles. If the label of $l$ is $y_{l_2}\dots
y_{l_m}$, then $s_1(D^{\prime})=x_{i_2}\cdot \, .\, .\, .\, \cdot
x_{i_m}$, $s_2(D^{\prime})=x_{l_2}\cdot \, .\, .\, .\, \cdot
x_{l_m}$ and $s_1(D^{\prime\prime})=x_{i_1}\cdot x_{l_2}\cdot \, .\,
.\, .\, \cdot x_{l_m}\cdot x_{i_{m+1}}\cdot \, .\, .\, .\, \cdot
x_{i_n}$, $s_2(D^{\prime\prime})=x_{j_1}\cdot \, .\, .\, .\, \cdot
x_{j_n}$ if $v_2\in \partial_l D$, and $s_1(D^{\prime})=x_{i_2}\cdot
\, .\, .\, .\, \cdot x_{i_n}$, $s_2(D^{\prime})=x_{l_2}\cdot \, .\,
.\, .\, \cdot x_{l_m}\cdot x_{j_{m+1}}\cdot \, .\, .\, .\, \cdot
x_{j_n}$ and $s_1(D^{\prime\prime})=x_{i_1}\cdot x_{l_2}\cdot \, .\,
.\, .\, \cdot x_{l_m}$, $s_2(D^{\prime\prime})=x_{j_1}\cdot \, .\,
.\, .\, \cdot x_{j_m}$ if $v_2\in \partial_r D$. If the end vertex
$v_m$ of $D^{\prime}$ (resp., $D^{\prime\prime}$) in Case $l$
(resp., in Case $r$) is not locally maximal, then we perform an
admissible transformation $II$ of $D^{\prime}$ (resp.,
$D^{\prime\prime}$) which increases neither the number of
quadrangles nor the number of locally maximal vertices.
\\

\begin{picture}(0,300)
\put(200,-15){$\text{Fig.}\, 6$} \put(116,28){$A$}
\put(70,210){\vector(1,1){20}}\put(90,230){\vector(1,1){20}}\put(120,260){\vector(1,1){20}}
\put(140,280){\circle*{3}} \put(136,285){$B$}
\put(90,70){\circle*{3}}\put(90,230){\circle*{3}}\put(90,70){\vector(0,1){30}}
\put(90,100){\vector(0,1){30}}\put(90,135){$\cdot$}
\put(90,140){$\cdot$}\put(90,145){$\cdot$}
\put(90,150){$\cdot$}\put(90,170){\vector(0,1){30}}
\put(90,200){\vector(0,1){30}}\put(90,155){$\cdot$}
\put(90,160){$\cdot$}\put(113,252){$\cdot$}\put(116,255){$\cdot$}
\put(110,249){$\cdot$}\put(52,192){$\cdot$}\put(57,197){$\cdot$}\put(62,202){$\cdot$}
\put(10,150){\vector(1,1){20}}\put(30,170){\vector(1,1){20}}
\put(10,150){\circle*{3}}\put(30,130){\vector(-1,1){20}}
\put(50,110){\vector(-1,1){20}}\put(90,70){\vector(-1,1){20}}
\put(55,100){$\cdot$}\put(60,95){$\cdot$}\put(65,90){$\cdot$}
\put(79,65){$v_1$}\put(120,40){\vector(-1,1){30}}
\put(115,272){$y_{i_n}$}\put(95,48){$y_{i_1}$}\put(92,80){$e_{2}$}
\put(92,110){$e_{3}$}\put(92,210){$e_{m}$}
\put(63,222){$y_{i_m}$}\put(100,233){$y_{i_{m+1}}$}\put(77,234){$v_{m}$}
\put(52,150){$D^{\prime}$}\put(112,150){$D^{\prime\prime}$}
\put(145,260){$y_{j_n}$}\put(130,43){$y_{j_1}$}\put(120,40){\vector(1,1){20}}
\put(120,40){\circle*{3}}\put(140,60){\circle*{3}}
\put(140,60){\vector(0,1){30}}\put(140,90){\vector(0,1){30}}
\put(140,130){$\cdot$}\put(140,135){$\cdot$} \put(140,140){$\cdot$}
\put(140,145){$\cdot$}\put(140,160){\vector(0,1){30}}\put(140,195){$\cdot$}
\put(140,200){$\cdot$}\put(140,205){$\cdot$}\put(140,210){$\cdot$}
\put(140,220){\vector(0,1){30}}\put(140,250){\vector(0,1){30}}\put(286,28){$A$}
\put(250,65){$v_1$}\put(245,90){$y_{i_2}$}\put(245,255){$y_{i_n}$}\put(268,47){$y_{i_1}$}
\put(278,80){$e_{2}$}\put(300,44){$y_{j_1}$}\put(280,280){$y_{j_n}$}
\put(353,179){$e_{m}$}\put(383,164){$y_{j_m}$}\put(383,209){$y_{j_{m+1}}$}
\put(383,187){$v_{m}$}\put(320,210){$D^{\prime}$}\put(320,80){$D^{\prime\prime}$}
\put(290,40){\circle*{3}}\put(290,40){\vector(-1,1){30}}\put(260,70){\circle*{3}}
\put(260,70){\vector(0,1){40}}\put(260,120){$\cdot$}\put(260,125){$\cdot$}
\put(260,130){$\cdot$}\put(260,115){$\cdot$}\put(260,140){\vector(0,1){30}}
\put(259,180){$\cdot$}\put(259,185){$\cdot$} \put(259,190){$\cdot$}
\put(259,195){$\cdot$}\put(259,200){$\cdot$}\put(260,215){\vector(0,1){30}}\put(260,245){\vector(0,1){30}}
\put(254,278){$B$}\put(260,275){\circle*{3}}\put(300,275){\vector(-1,0){40}}
\put(315,275){$\dots$}\put(380,275){\vector(-1,0){40}}\put(380,275){\circle*{3}}\put(380,225){\vector(0,1){50}}
\put(260,70){\vector(1,1){30}}\put(290,100){\vector(1,1){25}}
\put(330,140){\vector(1,1){20}}\put(325,133){$\cdot$}\put(322,130){$\cdot$}\put(319,127){$\cdot$}
\put(350,160){\vector(1,1){30}} \put(380,190){\vector(0,1){40}}
\put(380,150){\vector(0,1){40}}\put(380,190){\circle*{3}}
\put(380,40){\vector(0,1){40}}\put(380,80){\vector(0,1){40}}
\put(380,126){$\cdot$}\put(380,136){$\cdot$}\put(380,131){$\cdot$}
\put(290,40){\vector(1,0){30}}
\put(330,40){$\dots$}\put(350,40){\vector(1,0){30}}\put(380,40){\circle*{3}}
\put(100,5){$\text{Case}\, l$}\put(300,10){$\text{Case}\, r$}
\end{picture} \\ \\ \\

To complete the proof of Theorem \ref{main1}, we use two inductions:
the first one on the number $K$ of the invisible locally maximal
vertices of admissible disc diagrams satisfying conditions $(i)$ --
$(iii)$, and the second one on the number $k$ of quadrangles
entering into diagrams. As it was shown above, in each step of
inductions we can find a path $l$ dividing the admissible disc
diagram $D$ satisfying conditions $(i)$ -- $(iii)$ into two
subdiagrams $D^{\prime}$ and $D^{\prime\prime}$ such that either
(after an admissible transformation) $D^{\prime}$ and
$D^{\prime\prime}$ have not more than $K$ invisible locally maximal
vertices and have strictly less than $k$ quadrangles, or
$D^{\prime}$ and $D^{\prime\prime}$ have strictly less than $K$
invisible locally maximal vertices (if $K>0$). By inductive
assumptions, we have $s_1(D^{\prime})\sim s_2(D^{\prime})$ and
$s_1(D^{\prime\prime})\sim s_2(D^{\prime\prime})$. Therefore, by
Lemmas \ref{eqival2} and \ref{eqival3},
we have $s_1(D)\sim s_2(D)$. \qed \\

\subsection{Stability of the factorization semigroups over finite
$C$-groups} \label{st} Recall that the factorization semigroup
$S(G,O)$ over an equipped group $(G,O)$ is called {\it stable} if
there is an element $s\in S(G,O)$ such that $s_1\cdot s=s_2\cdot s$
for any two elements $s_1,s_2\in S(G,O)$ such that
$\alpha_G(s_1)=\alpha_G(s_2)$ and $\tau(s_1)=\tau(s_2)$.

Let $s_1,s_2$ be two equivalent elements of the factorization
semigroup $S(G,Y)$ of a $C$-group $(G,Y)$. Denote by $e(s_1,s_2)$
the smallest number $k$ such that there is an element $s\in S(G,Y)$
of length $k$ and such that $s_1\cdot s=s_2\cdot s$.

As above, let $\Gamma= \Gamma_1\sqcup\dots \sqcup \Gamma_m$ be the
decomposition of the finite $C$-graph of a $C$-group
$(G,Y)=G_{\Gamma}$ into the disjoint union of its connected
components, $V_i=\{ v_{i,1}, \dots , v_{i,n_i}\}$ be the set of
vertices of $\Gamma_i$, and  $p_i$ be the period of vertices
$v_{i,j}$. Any element $s\in S(G,Y)$ can be written in the form:
$\displaystyle s=\prod _{v_{i,j}\in \Gamma}x_{v_{i,j}}^{k_{i,j}}$.

Consider a set
$$R=\{  \displaystyle s=\prod
_{v_{i,j}\in \Gamma}x_{v_{i,j}}^{k_{i,j}}\in S(G,Y)\mid
\sum_{j=1}^{n_i}k_{i,j}\leq p_in_i,\, \, i=1,\dots ,m \},$$ where
$p_{i}$ are the periods of the $C$-generators $x_{v_{i,j}}$ of
$G_{\Gamma}$ and let $$E=\{ (s_1,s_2)\in R^2\mid s_1\sim s_2 \}.$$
Put
$$ \displaystyle e_{\Gamma}=\max_{(s_1,s_2)\in E}e(s_1,s_2).$$

\begin{lem} \label{divis1}
We have $s_1\cdot s_{\Gamma}^{e_{\Gamma}}=s_2\cdot
s_{\Gamma}^{e_{\Gamma}}$ for all $(s_1,s_2)\in E$, where
$s_{\Gamma}$ is the canonical element of $S(G,Y)$.
\end{lem}
\proof For any $(s_1,s_2)\in E$ there is an element $s_{1,2}$ of
$ln(s_{1,2})\leq e_{\Gamma}$ such that $s_1\cdot s_{1,2}=s_2\cdot
s_{1,2}$. By Lemma \ref{divis}, there is an element
$s_{1,2}^{\prime}$ such that $s_{1,2}\cdot
s_{1,2}^{\prime}=s_{\Gamma}^{e_{\Gamma}}$. Therefore
$$s_1\cdot s_{\Gamma}^{e_{\Gamma}}=s_1\cdot s_{1,2}\cdot
s_{1,2}^{\prime}=s_2\cdot s_{1,2}\cdot s_{1,2}^{\prime}=s_2\cdot
s_{\Gamma}^{e_{\Gamma}}.$$

\begin{thm} \label{stable}
Let $(G,Y)=G_{\Gamma}$ be a finite $C$-group. Then the semigroup
$S(G,Y)$ is stable and $s_{\Gamma}^{e_{\Gamma}}$ is a stabilizing
element.
\end{thm}
\proof

Let $s_1,s_2$ be two elements of the factorization semigroup
$S(G,Y)$ such that $\alpha_G(s_1)=\alpha_G(s_2)$ and
$\tau(s_1)=\tau(s_2)$. Write them in the form
$$s_1=\prod _{v_{i,j}\in
\Gamma}x_{v_{i,j}}^{k_{i,j,1}},\qquad s_2=\prod _{v_{i,j}\in
\Gamma}x_{v_{i,j}}^{k_{i,j,2}}.$$

Since $\tau(s_1)=\tau(s_2)$, we have
$$\sum_{j=1}^{n_i}k_{i,j,1}=\sum_{j=1}^{n_i}k_{i,j,2},\quad i=1,\dots, m.$$
If $\sum_{j=1}^{n_i}k_{i,j,1}\geq p_in_i$ for some $i$, then there
is $j_1$ (resp., $j_2$) such that $k_{i,j_1,1}\geq p_i$ (resp.,
$k_{i,j_2,2}\geq p_i$). If again
$$k_{i,j_1,1}-p_i+\sum_{j=1,j\neq j_1}^{n_i}k_{i,j,1}\geq p_in_i $$
(resp., $$k_{i,j_2,2}-p_i+\sum_{j=1,j\neq j_2}^{n_i}k_{i,j,2}\geq
p_in_i ),$$ then either there is $j^{\prime}_1$ (resp.,
$j^{\prime}_2$) such that $k_{i,j^{\prime}_1,1}\geq p_i$ (resp.,
$k_{i,j^{\prime}_2,2}\geq p_i$) or $k_{i,j_1,1}-p_i\geq p_i$ (resp.,
$k_{i,j_2,2}-p_i\geq p_i$). Continuing this process, as a result, we
obtain that the elements $s_1$ and $s_2$ can be written in the
following form:
$$s_1=\prod _{i=1}^{m}\prod _{j=1}^{n_i}x_{v_{i,j}}^{p_ic_{i,j,1}+r_{i,j,1}}=
(\prod _{i=1}^{m}\prod _{j=1}^{n_i}x_{v_{i,j}}^{p_ic_{i,j,1}})\cdot
(\prod _{i=1}^{m}\prod _{j=1}^{n_i}x_{v_{i,j}}^{r_{i,j,1}}),$$
$$s_2=\prod _{i=1}^{m}\prod _{j=1}^{n_i}x_{v_{i,j}}^{p_ic_{i,j,2}+r_{i,j,2}}=
(\prod _{i=1}^{m}\prod _{j=1}^{n_i}x_{v_{i,j}}^{p_ic_{i,j,2}})\cdot
(\prod _{i=1}^{m}\prod _{j=1}^{n_i}x_{v_{i,j}}^{r_{i,j,2}}), $$
where
$$\sum_{j=1}^{n_i}r_{i,j,1}=\sum_{j=1}^{n_i}r_{i,j,2}<n_ip_i,\quad i=1,\dots, m,$$
and
$$\sum_{j=1}^{n_i}c_{i,j,1}=\sum_{j=1}^{n_i}c_{i,j,2},\quad i=1,\dots, m.$$
Since
$\alpha_G(x_{v_{i,j_1}}^{p_i})=\alpha_G(x_{v_{i,j_2}}^{p_i})=y_{v_{i,1}}^{p_i}$
for $1\leq j_1\leq j_2\leq n_i$ and for all $i$ and since
$\alpha_G(s_1)=\alpha_G(s_2)$, we have
$$\alpha_G(\prod _{i=1}^{m}\prod _{j=1}^{n_i}x_{v_{i,j}}^{r_{i,j,1}})
=\alpha_G(\prod _{i=1}^{m}\prod
_{j=1}^{n_i}x_{v_{i,j}}^{r_{i,j,2}}).$$ Then, by Theorem
\ref{main1}, we have
$$\prod _{i=1}^{m}\prod _{j=1}^{n_i}x_{v_{i,j}}^{r_{i,j,1}}\sim
\prod _{i=1}^{m}\prod _{j=1}^{n_i}x_{v_{i,j}}^{r_{i,j,2}},$$
and by Lemma \ref{divis1},
$$(\prod _{i=1}^{m}\prod _{j=1}^{n_i}x_{v_{i,j}}^{r_{i,j,1}})\cdot s_{\Gamma}^{e_{\Gamma}}
= (\prod _{i=1}^{m}\prod _{j=1}^{n_i}x_{v_{i,j}}^{r_{i,j,2}})\cdot
s_{\Gamma}^{e_{\Gamma}}.$$

Note that $s_1 \cdot  s_{\Gamma}^{e_{\Gamma}}$ (resp., $s_2 \cdot
s_{\Gamma}^{e_{\Gamma}}$) belongs to $S(G,Y)^G$. Therefore, by Lemma
\ref{fried},
$$s_1 \cdot  s_{\Gamma}^{e_{\Gamma}}=
(\prod _{i=1}^{m}\prod _{j=1}^{n_i}x_{v_{i,j}}^{p_ic_{i,j,1}})\cdot
(\prod _{i=1}^{m}\prod _{j=1}^{n_i}x_{v_{i,j}}^{r_{i,j,1}}) \cdot
s_{\Gamma}^{e_{\Gamma}}=(\prod
_{i=1}^{m}x_{v_{i,1}}^{p_ic_{i}})\cdot (\prod _{i=1}^{m}\prod
_{j=1}^{n_i}x_{v_{i,j}}^{r_{i,j,1}}) \cdot s_{\Gamma}^{e_{\Gamma}}$$
and
$$s_2 \cdot  s_{\Gamma}^{e_{\Gamma}}=
(\prod _{i=1}^{m}\prod _{j=1}^{n_i}x_{v_{i,j}}^{p_ic_{i,j,2}})\cdot
(\prod _{i=1}^{m}\prod _{j=1}^{n_i}x_{v_{i,j}}^{r_{i,j,2}}) \cdot
s_{\Gamma}^{e_{\Gamma}}=(\prod
_{i=1}^{m}x_{v_{i,1}}^{p_ic_{i}})\cdot (\prod _{i=1}^{m}\prod
_{j=1}^{n_i}x_{v_{i,j}}^{r_{i,j,2}}) \cdot
s_{\Gamma}^{e_{\Gamma}},$$ where
$c_{i}=\sum_{j=1}^{n_i}c_{i,j,1}=\sum_{j=1}^{n_i}c_{i,j,2}$, and
hence $s_1 \cdot  s_{\Gamma}^{e_{\Gamma}}=s_2 \cdot
s_{\Gamma}^{e_{\Gamma}}$. \qed \\
\begin{thm} \label{stabel2}
Let $(G,O)$ be an equipped finite group. The semigroup $S(G,O)$ {\rm
(}resp., $S(G,O)^G$, $S(G,O)_{\bf{1}}$, and  $S(G,O)^G_{\bf{1}}$
{\rm )} over the group $G$ is stable if and only if the ambiguity
index $a_{(G,O)}=1$.
\end{thm}
\proof The equipped group $(G,O)$ is equivalent to the finite
$C$-group $(G_{\Gamma},Y)$, where $\Gamma=\Gamma_{(G,O)}$, and there
is an epimorphism $\beta_{(G,O)}:G_{\Gamma}\to G$, such that $H=\ker
\beta\subset Z(G_{\Gamma})$. By Claim \ref{iso3}, the semigroups
$S(G,O)$ and $S(G_{\Gamma},Y)$ are naturally isomorphic. Let
$\alpha_{G_{\Gamma}}$ and $\alpha_G$ be respectively the product
homomorphisms of $S(G,O)$ to $G_{\Gamma}$ and $G$. Note that, by
Theorem \ref{stable}, for any positive integer $c$ the element
$s_{\Gamma}^{ce_{\Gamma}}$ is a stabilizing element of
$S(G_{\Gamma},Y)$ over the group $G_{\Gamma}$ and there is a
positive integer $c_0$ such that
$\alpha_G(s_{\Gamma}^{c_0e_{\Gamma}})={\bf{1}}\in G$, since $G$ is a
finite group.

If $a_{(G,O)}=|H\cap [G_{\Gamma},G_{\Gamma}]|=1$, then for
$s_1,s_2\in S(G,O)$ such that $\tau(s_1)=\tau(s_2)$, we have
$\alpha_G(s_1)=\alpha_G(s_2)$ if and only if
$\alpha_{G_{\Gamma}}(s_1)=\alpha_{G_{\Gamma}}(s_2)$. Therefore if
$a_{(G,O)}=1$, then the semigroups $S(G,O)$, $S(G,O)^G$,
$S(G,O)_{\bf{1}}$, and  $S(G,O)^G_{\bf{1}}$ over the group $G$ are
stable and $s_{\Gamma}^{c_0e_{\Gamma}}$ is one of their stabilizing
elements.

If $a_{(G,O)}>1$, then there is an element $g\in H\cap
[G_{\Gamma},G_{\Gamma}]$ such that $g\neq {\bf{1}}$. By Lemma
\ref{simple2}, there are two elements $s_1, s_2\in S(G_{\Gamma},Y)$
such that $g=\alpha_{G_{\Gamma}}(s_1)\alpha_{G_{\Gamma}}(s_2)^{-1}$.
Therefore, $s_1\not\sim s_2$, but $\alpha_{G}(s_1)=\alpha_{G}(s_2)$
and $\tau(s_1)=\tau(s_2)$, that is, the semigroups $S(G,O)$,
$S(G,O)^G$, $S(G,O)_{\bf{1}}$, and  $S(G,O)^G_{\bf{1}}$ over the
group $G$ are not stable  (multiplying $s_1$ and $s_2$ by some
element $s$, we can
assume that $s_1,s_2\in S(G,O)_{\bf{1}}^G$ ). \qed \\

\section{Uniqueness of factorizations in the case of big enough number of factors }
\subsection{The case of finite $C$-groups}
 Let $G_{\Gamma}=(G,Y)$ be a finite
$C$-group and $m$ the number of connected components of the
$C$-graph $\Gamma$.
\begin{thm} \label{Cuniq} For each finite $C$-group $G_{\Gamma}=(G,Y)$, there is a constant $T\in \mathbb N$ such that if
elements $s_1,s_2\in S(G,Y)^G$ satisfy the following conditions:
\begin{itemize}
\item[$(i)$] $\tau_i(s_1)\geq T$ for $i=1,\dots, m$;
\item[$(ii)$] $\alpha_G(s_1)=\alpha_G(s_2)$, \end{itemize} then $s_1=s_2$.
\end{thm}
\proof Put $\alpha=\alpha_G$ and note, first of all, that if
$\alpha(s_1)=\alpha(s_2)$, then $\tau_i(s_1)=\tau_i(s_2)$ for
$i=1,\dots, m$, since $\tau_i(s_j)=ab_i(\alpha(s_j))$.

Let us denote by $d_{\Gamma}$ the diameter of $\Gamma$ and, using
the notations of section \ref{st}, denote
$$\displaystyle T_1=(d_{\Gamma}+1)\max_{1\leq i\leq m}n_ip_i +1.$$

Let us show that any integer $T\geq T_1$ satisfies the conditions of
Theorem \ref{Cuniq}. Indeed let $\tau_i(s_1)\geq T_1$ for
$i=1,\dots, m$. Then, by Corollary \ref{can}, the elements $s_1$ and
$s_2$ can be written in the form: $s_1=s_{\Gamma}^{d_{\Gamma}}\cdot
s_1^{\prime}$ and $s_2=s_{\Gamma}^{d_{\Gamma}}\cdot s_2^{\prime}$.
We have
$$\alpha(s_1^{\prime})=\alpha(s_1)\alpha(s_{\Gamma}^{d_{\Gamma}})^{-1}=
\alpha(s_2)\alpha(s_{\Gamma}^{d_{\Gamma}})^{-1}=\alpha(s_2^{\prime})$$
Therefore $s_1^{\prime}\sim s_2^{\prime}$ by Theorem \ref{main1} and
hence $s_1=s_2$ by Theorem \ref{stable}. \qed \\

\subsection{The case of equipped finite groups} \label{TT}
Let $(G,O)$ be an equipped group equivalent to a finite $C$-group
$G_{\Gamma}$, $\beta_{(G,O)}:G_{\Gamma}\to (G,O)$ the natural
epimorphism of equipped groups, $H=\ker \beta_{(G,O)}$,
$a_{(G,O)}=|H\cap [G_{\Gamma},G_{\Gamma}]|$ the ambiguity index of
$(G,O)$, and $m$ the number of connected components of $\Gamma$.

\begin{thm} \label{Cuniq1} For each equipped finite group $(G,O)$ there is a constant $T=T_{(G,O)}$
such that if for an element $s_1\in S(G,O)^G$ the $i$th type
$\tau_i(s_1)\geq T$ for all $i=1,\dots, m$, then there are
$a_{(G,O)}$ elements  $s_1,\dots, s_{a_{(G,O)}}\in S(G,O)^G$ such
that
\begin{itemize}
\item[$(i)$] $s_i\neq s_j$ for $1\leq i< j\leq a_{(G,O)}$;
\item[$(ii)$] $\tau(s_i)=\tau(s_1)$ for $1\leq i\leq a_{(G,O)}$;
\item[$(iii)$] $\alpha_G(s_i)=\alpha_G(s_1)$ for $1\leq i\leq a_{(G,O)}$;
\item[$(iv)$] if $s\in S(G,O)^G$ is such that $\tau(s)=\tau(s_1)$ and $\alpha_G(s)=\alpha_G(s_1)$, then
$s=s_i$ for some $i$, $1\leq i\leq  a_{(G,O)}$.
\end{itemize}
\end{thm}
\proof By Lemma \ref{simple2}, for $n\geq z_H$ each element $g_i\in
H\cap [G_{\Gamma},G_{\Gamma}]$, (where $g_1={\bf{1}}$) can be
presented in the form
$g_i=\alpha_{G_{\Gamma}}(s_{i,n})\alpha_{G_{\Gamma}}(s_{\Gamma}^n)^{-1}$,
where $s_{i,n}\in S(G_{\Gamma},Y)=S(G,O)$ (in particular,
$s_{1,n}=s_{\Gamma}^n$). We have $s_{i,n}\neq s_{j,n}$ for $1\leq i<
j\leq a_{(G,O)}$, $\tau(s_{i,n})=\tau(s_{1,n})$ and
$\alpha_G(s_{i,n})=\alpha_G(s_{1,n})$ for $1\leq i\leq a_{(G,O)}$.

Put $T_2=(\max(d_{\Gamma},z_H)+1)\max_{1\leq i\leq m}n_ip_i +1$. By
Corollary \ref{can}, if $T\geq T_2$, then the element $s_1$ can be
written in the form: $s_1=s_{\Gamma}^{z_{H}}\cdot s_1^{\prime}$.
Denote $s_i:=s_{i,z_H}^{z_{H}}\cdot s_1^{\prime}$. It is easy to see
that the elements $s_1,\dots, s_{a_{(G,O)}}\in S(G,O)^G$ satisfy
conditions $(i)$ -- $(iii)$ and if $s\in S(G,O)^G$ is such that
$\tau(s)=\tau(s_1)$ and $\alpha_G(s)=\alpha_G(s_1)$, then
$\alpha_{\Gamma}(s)\alpha_G(s_1)^{-1}\in H\cap
[G_{\Gamma},G_{\Gamma}]$ and hence, by Theorem \ref{Cuniq}, $s=s_i$
for some $i$, $1\leq i\leq  a_{(G,O)}$. \qed \\

Let $(G,O)$ be an equipped group, $O=C_1\sqcup\dots\sqcup C_m$,
where $C_i$ are conjugcy classes of $G$. Let for some $k<m$ the
elements of the set $O^{\prime}=C_1\sqcup\dots\sqcup C_k$ generate
the group $G$. The embedding $i:O^{\prime}\hookrightarrow O$ defines
subgraphs
$\widetilde{\Gamma}^{\prime}\subset\Gamma^{\prime}=\Gamma_1\sqcup
\dots \sqcup \Gamma_k$ of the $C$-graph $\Gamma_{(G,O)}=\Gamma$,
where $\widetilde{\Gamma}^{\prime}$ is the $C$-graph of the equipped
group $(G,O^{\prime})$. The subgraph $\Gamma^{\prime}$ is ample,
 since the the elements of  $O^{\prime}$ generate the
group $G$.

Let $G_{\Gamma}$ and $G_{\widetilde{\Gamma}^{\prime}}$ be $C$-groups
equivalent respectively to $(G,O)$ and $(G,O^{\prime})$, and
$S(G_{\Gamma},Y)$ and
$S(G_{\widetilde{\Gamma}^{\prime}},Y^{\prime})$ their factorization
semigroups. The embedding $i:O^{\prime}\hookrightarrow O$ defines a
homomorphism $i_*:G_{\widetilde{\Gamma}^{\prime}}\to  G_{\Gamma}$ of
$C$-groups and an embedding $$i_*:S(G,O^{\prime})\simeq
S(G_{\widetilde{\Gamma}^{\prime}},Y^{\prime}) \hookrightarrow
S(G_{\Gamma},Y)\simeq S(G,O)$$ of semigroups such that
$i_*(\alpha_{G_{\widetilde{\Gamma}^{\prime}}}(s))=\alpha_{G_{\Gamma}}(i_*(s))$
for all $s\in S(G_{\widetilde{\Gamma}^{\prime}},Y^{\prime})$.

\begin{thm} \label{Cuniq2} Let $G$ be a finite group and
$O^{\prime}\subset O$ be two its equipments such that the elements
of $O^{\prime}$ generate the group $G$. Then, in notations used
above, there is a constant $T=T_{(O,O^{\prime})}$ such that if for
an element $s_1\in S(G,O)^G$ the $i$th type $\tau_i(s_1)\geq T$ for
all $i=1,\dots, k$, then there are not more than
$a_{(G,O^{\prime})}$ elements $s_1,\dots, s_{n}\in S(G,O)^G$ such
that
\begin{itemize}
\item[$(i)$] $s_i\neq s_j$ for $1\leq i< j\leq n$;
\item[$(ii)$] $\tau(s_i)=\tau(s_1)$ for $1\leq i\leq n$;
\item[$(iii)$] $\alpha_G(s_i)=\alpha_G(s_1)$ for $1\leq i\leq n$,
\end{itemize}
where $a_{(G,O^{\prime})}$ is the ambiguity index of
$(G,O^{\prime})$.
\end{thm}
\proof Let $\displaystyle T_1=2d_{\Gamma^{\prime}}\max_{1\leq i\leq
k} n_ip_i$, where $d_{\Gamma^{\prime}}$ is the diameter of $\Gamma$
with respect to the ample subgraph $\Gamma^{\prime}$, $p_i$ and
$n_i$, respectively, be the period and the number of vertices of the
connected component $\Gamma_i$ of the graph $\Gamma$, and $T_2$ be a
constant the existence of which for the equipped group
$(G,O^{\prime})$ is claimed in Theorem \ref{Cuniq1}. Put
$T=\max(T_1,T_2)$. By Proposition \ref{rep2}, if for $s\in S(G,Y)^G$
its $i$th type $\tau_i(s)\geq T$ for all $i\leq k$, then $s$ can be
written in the form: $s=(x_{v_{{k+1},1}}^{a_{k+1}}\cdot \, .\, .\,
.\, \cdot x_{v_{m,1}}^{a_m})\cdot s_{\Gamma^{\prime}}^d\cdot s_1$,
where $a_i=\tau_i(s)$ for $i=k+1,\dots, m$ and $s_1\in
S(G,O^{\prime})$. Now, Theorem \ref{Cuniq2} follows from  Theorem
\ref{Cuniq1} applied to the element $s_{\Gamma^{\prime}}^d\cdot
s_1\in S(G,O^{\prime})^G$. \qed

\subsection{Generating functions} For each $\overline k=(k_1,\dots,k_m)\in
\mathbb Z_{\geq 0}^m$ denote by $h_{\overline k}$ the number of
elements $s\in S(G,O)^G_{\bf{1}}$ with $\tau(s)=\overline k$ and
associate to $(G,O)$ a power series $\displaystyle
\chi_{(G,O)}(t_1,\dots,t_m):=\sum_{\overline k\in \mathbb Z_{\geq
0}^m}h_{\overline k}t_1^{k_1}\dots t_m^{k_m}.$

\begin{prop}\label{xyz}
Let $G$ be a finite group and $C$ its conjugacy class such that the
elements of $C$ generate the group $G$. Then $\chi_{(G,C)}(t)$ is a
rational function.
\end{prop}
\proof Let $p$ be the period and $n$ the number of the vertices of
$\Gamma_{(G,C)}$, and $a=a_{(G,C))}$ index of ambiguity of $(G,C)$.
Consider a set of integers $$R=\{ r\mid 0\leq r\leq pn-1\, \,
\text{and}\, \, \exists \, i\equiv r(\text{mod}\, np)\, \,
\text{such that}\, \, h_i>0\}$$ and let $i_r$ be the smallest $i$
for which $i\equiv r(\text{mod}\, np)$ and $h_i>0$. Next, for each
$r\in R$ choose a representative $s_r\in S(G,C)_{\bf{1}}^G$ of
length $ln(s_r)=i_r$ and choose a constant $M$ such that
$\displaystyle T=2pnM>\max_{r\in R}i_r$ is a constant the existence
of which is proved in Theorem \ref{Cuniq1}.

Write $\chi_{(G,C)}(t)$ in the form:
$\chi_{(G,C)}(t)=\chi_{<T}(t)+\chi_{\geq T}(t)$, where
$\chi_{<T}(t)=\sum_{i=1}^{T-1}h_it^i$. Note that $\chi_{<T}(t)$ is a
polynomial. The function $\chi_{\geq T}(t)$ can be written in the
form $$\chi_{\geq T}(t)=\sum_{r\in
R}\sum_{j=0}^{\infty}h_{r+T+jpn}t^{r+T+jpn}.
$$
By  Theorem \ref{Cuniq1}, we have $h_{r+T+jpn}=a$ for each $r\in R$
and each $j\geq 0$, therefore
$$\chi_{\geq T}(t)=at^T\sum_{r\in R}t^{r}\sum_{j=0}^{\infty}t^{jpn}=
\frac{ at^T}{1-t^{pn}}\sum_{r\in R}t^{r}.\qquad \qed$$

\begin{ex} {\rm If $G=\Sigma_n$ is the symmetric group and $O$ the set of
transpositions, then by Clebsch -- Hurwitz Theorem},
$\chi_{(G,O)}=\frac{t^{2(n-1)}}{1-t^2}$. \end{ex}

A generalization of Proposition \ref{xyz} is the following
\begin{thm}
Let $G$ be a finite group and $O=C_1\sqcup\dots \sqcup C_m$ a
disjoint union of its conjugacy classes such that the elements of
each class $C_i$ generate the group $G$. Then
$\chi_{(G,O)}(t_1,\dots,t_m)$ is a rational function.
\end{thm}
\proof It is similar to the proof of Proposition \ref{xyz}. We must
note only that if $k_i$ is big enough, then it follows from  Lemma
\ref{x1} and Theorem \ref{Cuniq2} that
$$0\leq h_{(k_1,\dots,k_{i-1},k_i+p_in_i,k_{i+1}, \dots, k_m)}\leq
h_{(k_1,\dots,k_{i-1},k_i,k_{i+1}, \dots, k_m)}. \qquad \qed$$

\begin{que}  Is $\chi_{(G,O)}(\overline t)$ a rational function for
any equipped finite group $(G,O)${\rm ?}
\end{que}

\section{Computation of the ambiguity index}
\subsection{The word problem for finite $C$-groups} In this subsection we prove
\begin{thm} \label{algol} The finite $C$-groups have solvable word problem. \end{thm}
\proof Let $\Gamma=\Gamma_1\sqcup\dots\sqcup \Gamma_m$ be the
$C$-graph of a finite $C$-group $(G_{\Gamma},Y)$. By Claim
\ref{ppp}, without loss of generality, we can assume that for each
connected component $\Gamma_i$ of  $\Gamma$ the period $p_i$ of its
vertices is greater than one.

By Proposition \ref{fincom}, applying the relations of $C$-group
given by $C$-graph $\Gamma$, any word in letters of $Y$  can be
transformed into a word of the following form:
\begin{equation}\label{ccom1} w=c^{-1}y_{1,1}^{t_{1}}y_{2,1}^{t_{2}}\dots y_{m,1}^{t_{m}}
\prod_{i=1}^{m}y_{i,1}^{k_ip_i}y_{i,1}^{a_{i,1}}y_{i,2}^{a_{i,2}}\dots
y_{i,n_i}^{a_{i,n_i}},
\end{equation}
 where $c$ is the canonical element of $G_{\Gamma}$ and the integers $k_i$ and $a_{i,j}$
satisfy the following relations and inequalities
\begin{equation} \label{xxxx1}
\sum_{j=1}^{n_i}a_{i,j}+k_ip_i=n_ip_i,\end{equation}
\begin{equation}\label{yyyy1}
 0<a_{i,j}\leq p_i-1,  \qquad 0\leq k_i< n_i.
\end{equation}

If two words $w_1$ and $w_2$  represent the same element $g\in
G_{\Gamma}$, then $ab(w_1)=ab(w_2)=(t_1,\dots,t_m)$. Therefore, to
prove Theorem, it suffices to show that there is a finite algorithm
solving the following problem:
\newline {\it To recognize when two words
\begin{equation}\label{ccom2} w_l=
\prod_{i=1}^{m}y_{i,1}^{a_{i,1,l}}y_{i,2}^{a_{i,2,l}}\dots
y_{i,n_i}^{a_{i,n_i,l}},\qquad l=1,2,
\end{equation}
satisfying the following relations and inequalities
\begin{equation} \label{xxxx2}
\sum_{j=1}^{n_i}a_{i,j,l}=k_{i,l}p_i,\end{equation}
\begin{equation}\label{yyyy2}
 0\leq a_{i,j,l}\leq p_i-1,  \qquad 0\leq k_{i,l}< n_i
\end{equation}
represent the same element in} $G_{\Gamma}$.

Note also that the words
\begin{equation}\label{commuta} w_l\prod_{i=1}^my_{i,1}^{-k_{i,l}p_i}, \qquad l=1,2,\end{equation}
where $w_l$ are words from (\ref{ccom2}), represent elements of
$[G_{\Gamma},G_{\Gamma}]$.

In the sequel, to simplify notations, we will consider only the case
when $m=1$, since the general case is similar. We will use the
following notations: $\Gamma=\Gamma_1$, $n:=n_1$, $p:=p_1$,
$\{v_1,\dots,v_n\}$ is the set of vertices of $\Gamma$ and $y_i$ is
the $C$-generator of $G_{\Gamma}$ corresponding to a vertex $v_i$.
Denote by $\mathcal R_{\Gamma}$ the set of relations in $G_{\Gamma}$
defined by the $C$-graph $\Gamma$ (see subsection \ref{C-rel}).
Remind that for each pair $(l_1,l_2)$ such that $1\leq l_1,l_2\leq
n$, there is unique $l_3$ (depending on $(l_1,l_2)$), $1\leq l_3\leq
n$, such that
\begin{equation} \label{rem}
y_{l_1}y_{l_2}y_{l_1}^{-1}y_{l_3}^{-1}\in \mathcal
R_{\Gamma}.\end{equation}

 Consider an equipped group $(\overline G=G_{\Gamma}/H,O)$, where $H$ is the subgroup of the center
$Z(G_{\Gamma})$ generated by the element $y_{1}^{p}=y_i^p$,
$i=1,\dots,n$. Denote by $f=f_{(G,O)}:G_{\Gamma}\to \overline G$ the
natural epimorphism and by $f^{\prime}=f_{\mid
[G_{\Gamma},G_{\Gamma}]}$ its restriction to
$[G_{\Gamma},G_{\Gamma}]$. The group $\overline G$ has the following
presentation:
$$\overline G=\langle y_1,\dots,y_n\mid R(\overline y)=1 \, \,
\text{for}\, \, R\in \mathcal R_{\Gamma}\, \, \text{and}\, \,
y_1^p=\dots =y_n^p=1\rangle .$$

By Lemma \ref{alg1}, we have the commutative diagram
\\

\begin{picture}(0,90)
\put(85,75){$1\longrightarrow  [G_{\Gamma},G_{\Gamma}]
\longrightarrow G \buildrel{ab}\over\longrightarrow  \mathbb Z
\longrightarrow 0$} \put(92,25){$1\longrightarrow [\overline
G,\overline  G] \longrightarrow \overline G
\buildrel{ab}\over\longrightarrow \mathbb Z/p\mathbb Z \to 0$}
\put(135,70){\vector(0,-1){30}} \put(126,55){$f^{\prime}$}
\put(187,70){\vector(0,-1){30}}
\put(189,55){$f$}\put(213,55){$f^{\prime\prime}$}
\put(225,70){\vector(0,-1){30}} \put(137,55){$\simeq$}
\put(390,55){$(*)$}
\end{picture}\newline
in which $f^{\prime}$ is an isomorphism, and $f$ and
$f^{\prime\prime}$ are epimorphisms. It follows from the above
considerations and diagram ($*$) that each element $g\in \overline
G$ can be represented by a word of the form:
\begin{equation}\label{ccom3} w=
y_{1}^{a_{1}}y_{2}^{a_{2}}\dots y_{n}^{a_{n}},
\end{equation}
where each integer $a_j$ satisfies the inequality
\begin{equation}\label{yyyy3}
 0\leq a_{j}\leq p-1
\end{equation}
and hence
\begin{equation} \label{xxxx3}
\sum_{j=1}^{n}a_{j}\leq n(p-1).\end{equation}

Denote by $W_{i}$ the set of positive words $w$ in letters of $Y$
whose lengths $ln(w)=i$ (in particular, $W_0$ consists of the empty
word), $W_{\leq k}=\bigcup_{i=0}^kW_i$, and by
$$\widetilde W_{\leq k}=\{
w\in W_{\leq k}\mid \text{each letter $y_j\in Y$ enters in $w$ less
than $p$ times}\}.$$ Therefore, to prove Theorem \ref{algol}, it
suffices to show that there is a finite algorithm recognizing when
two words $w_1,w_2\in \widetilde W_{\leq n(p-1)}$ represent the same
element of $\overline G$. Prior to describe a such algorithm, we
give two definitions. Let $w$ be a positive word in letters of $Y$
in which some letter $y_j\in Y$ enters sequentially $p$ times, that
is, $w=w^{\prime} y_j^pw^{\prime\prime}$. Then we call the word
$\overline w=w^{\prime}w^{\prime\prime}$  a {\it reduction} of $w$.
For each subset $V\subset W_{\leq k}$, denote by $\overline V=\{
\widetilde w\in \widetilde W_{\leq n(p-1)}\mid \exists w\in V\, \,
\text{such that}\, \, \widetilde w=\overline w\}.$

Let $\sigma_1,\dots, \sigma_{k-1}$ be a set of Artin generators of
the braid group $B_k$. Define an action of the group $B_k$  on the
set $W_i$ as follows: if $i\neq k$, then this action is trivial; and
if $i=k$, then for a word $w=y_{i_1}\dots y_{i_{j-1}}
y_{i_j}y_{i_{j+1}}y_{i_{j+2}}\dots y_{i_k}$ its image
$\sigma_j(w)=y_{i_1}\dots y_{i_{j-1}} y_{
i^{\prime}}y_{i_{j}}y_{i_{j+2}}\dots y_{i_k}$, where $ i^{\prime}$
depends on the pair $(i_j,i_{j+1})$ and it is defined by relation
(\ref{rem}) if we put $l_1=i_j$, $l_2=i_{j+1}$, and
$i^{\prime}=l_3$. Denote by $$B=B_0\times B_1\times \dots \times
B_{2n(p-1)},$$  where $B_0$ and $B_1$ are trivial groups. The
action, defined above, makes it possible to define an action of $B$
on $W_{\leq 2n(p-1)}=\bigcup_{i=0}^{2n(p-1)}W_i$. For each subset
$V$ of $W_{\leq 2n(p-1)}$ denote by $BV$ the union of the orbits of
the elements of $V$ under the action of $B$.

Let $N_1$ be the number of words in $\widetilde W_{\leq n(p-1)}$ and
let us numerate them, $w_1,\dots,w_{N_1}$, and consider each of them
as a subset of $\widetilde W_{\leq n(p-1)}$ (denote by $\widetilde
V_{j,1}=\{ w_j\}$). The algorithm can be described as follows. Let
in the end of $(k-1)$th step we have obtained a presentation of
$\widetilde W_{\leq n(p-1)}$ as a disjoint union of its subsets,
$$\widetilde W_{\leq n(p-1)}=\widetilde V_{1,k}\sqcup\dots\sqcup \widetilde V_{N_k,k}.$$
On the $k$th step for each pair $(i,j)\in \{ 1,\dots,N_k\}^2$, we
form subsets
$$\widetilde V_{i,k}\widetilde V_{j,k}=\{
w=w^{\prime}w^{\prime\prime}\in W_{\leq 2n(p-1)}\mid w^{\prime}\in
\widetilde V_{i,k},\, \,  w^{\prime\prime}\in \widetilde
V_{j,k}\}\subset W_{\leq 2n(p-1)}$$ and form the set of their orbits
$\{ B(\widetilde V_{i,k}\widetilde V_{j,k})\}$. Define an
equivalence relation induced by the following equivalence: two
orbits $B(\widetilde V_{i_1,k}\widetilde V_{j,k})$ and $B(\widetilde
V_{i_2,k}\widetilde V_{j,k})$ (resp., $B(\widetilde
V_{i,k}\widetilde V_{j_1,k})$ and $B(\widetilde V_{i,k}\widetilde
V_{j_2,k})$) are equivalent if there is $j_0$ (resp., $i_0$) such
that $$B(\widetilde V_{i_1,k}\widetilde V_{j_0,k})\cap B(\widetilde
V_{i_2,k}\widetilde V_{j_0,k})\neq \emptyset$$ (resp., $B(\widetilde
V_{i_0,k}\widetilde V_{j_1,k})\cap B(\widetilde V_{i_0,k}\widetilde
V_{j_2,k})\neq \emptyset$) and for each equivalence class we unite
the subsets $B(\widetilde V_{i,k}\widetilde V_{j,k})$ belonging to
this class. Denote the obtained subsets by $V_{1,k},\dots,
V_{N_{k+1},k}$ and put $\widetilde V_{i,k+1}:= \overline V_{i,k}$.

It is easy to see that $N_{k+1}\leq N_k$ and $\widetilde
V_{1,k+1}\sqcup\dots\sqcup \widetilde V_{N_{k+1},k+1}=\widetilde
W_{\leq n(p-1)}$. The algorithm is stopped if $N_{k+1}= N_k$ and two
words $w_l= y_{1}^{a_{1,l}}y_{2}^{a_{2,l}}\dots y_{n}^{a_{n,l}}$,
$l=1,2$, satisfying relation (\ref{xxxx2}) and inequality
(\ref{yyyy2}) (in the case $m=1$) represent the same element in
$G_{\Gamma}$ if and only if $w_1$ and $w_2$ belong to the same
subset $\widetilde V_i:=\widetilde V_{i,k+1}$ for some $i$. To show
this, let us introduce a group structure on the set $\widetilde G=\{
\widetilde V_i\}$. By definition, the product of $\widetilde V_i$
and $\widetilde V_j$ is $\widetilde V_i\widetilde
V_j=\overline{B(\widetilde V_i\widetilde V_j)}$. This product is
well defined by construction of subsets $\widetilde V_i$ and since
$N_{k+1}=N_{k}$. The unite in $\widetilde G$ is a subset $\widetilde
V_{i_0}$ containing the empty word, the inverse element of
$\widetilde V_{i}$ containing a word
$w=y_{1}^{a_{1}}y_{2}^{a_{2}}\dots y_{n}^{a_{n}}$ is the subset
$\widetilde V_{j}$ containing the reduction of the word
$w=y_{n}^{p-a_{n}}y_{n-1}^{p-a_{n-1}}\dots y_{1}^{p-a_{1}}$. Let us
renumber the subsets $\widetilde V_j$ so that for $i=1,\dots,n$ the
subset $\widetilde V_i$ contains the word $y_i$ and put $\widetilde
O=\{ \widetilde V_1,\dots, \widetilde V_n\}$. Then it is easy to see
that $\widetilde O$ is invariant under the inner automorphisms of
$\widetilde G$ and the $C$-graph of the equipped group $(\widetilde
G,\widetilde O)$ coincides with $\Gamma$. Moreover, since to
construct the elements of $\widetilde G$, we used only the relations
from $\mathcal R_{\Gamma}$ and relations $y_1^p=\dots =y_n^p=1$,
there is an epimorphism $h:(\widetilde G,\widetilde O)\to (\overline
G,\overline O)$ such that  the following  diagram \\ \\
\begin{picture}(0,75)
\put(155,52){$G_{\Gamma}$}
\put(170,55){\vector(1,0){50}}\put(225,52){$\widetilde G$}
\put(222,30){$h$} \put(167,48){\vector(1,-1){25}}
\put(227,48){\vector(-1,-1){25}}
\put(170,30){$f$}\put(192,10){$\overline
G$}\put(180,61){$f_{(\widetilde G,\widetilde O)}$}
\end{picture}
\newline  is commutative. Therefore $h([\widetilde G,\widetilde G])=[\overline G,\overline G]$ and since
$f_{(\widetilde G,\widetilde
O)}([G_{\Gamma},G_{\Gamma}])=[\widetilde G,\widetilde G]$ and
$f^{\prime}: [G_{\Gamma},G_{\Gamma}]\to [\overline G,\overline G]$
is an isomorphism, then $f_{(\widetilde G,\widetilde
O)}:[G_{\Gamma},G_{\Gamma}]\to [\widetilde G,\widetilde G]$ is also
an isomorphism. \qed
\subsection{Computation of the ambiguity index}
We say that an equipped finite group $(G,O)$ is {\it defined
efficiently} if there is a finite algorithm to enumerate the
elements of $G$ (for example, if $G$ is given by its Cayley graph,
or generators of $G$  as a subgroup of some symmetric group are
given) and a representative of each class $C_i\subset O$ is also
given.
\begin{prop} If  an equipped finite group $(G,O)$ is defined
efficiently, then there is a finite algorithm to compute the
ambiguity index $a_{(G,O)}$.
\end{prop}
\proof Since $(G,O)$ is given efficiently, then, obviously, there is
a finite algorithm to define completely the $C$-graph
$\Gamma=\Gamma_{(G,O)}$.

Using the description of $\Gamma$, as in the proof of Theorem
\ref{algol}, we enumerate the subsets $\widetilde V_i$ of the set of
words of the form (\ref{ccom2}) representing the same element in
$G_{\Gamma}$ and which satisfy relations (\ref{xxxx2}) and
inequalities (\ref{yyyy2}). After that for each $i$ we choose a word
$w_i\in \widetilde V_i$ and check if $w_i$ represents the unity in
$G$. The number of the sets $\widetilde V_i$ representing the unity
of $G$ is equal to $a_{(G,O)}$. \qed

 \ifx\undefined\bysame
\newcommand{\bysame}{\leavevmode\hbox to3em{\hrulefill}\,}
\fi

\end{document}